\documentclass[10pt]{article}
\usepackage{graphicx,amsthm,enumerate}
\usepackage{verbatim,amsmath,amscd,amssymb,color}
\usepackage{amsmath,amssymb,amsthm,color,eepic,eclarith}
\theoremstyle{plain}

\newcommand{\MECry}{\fsquare(0.3cm,i_1) \otimes \fsquare(0.3cm,i_2) \otimes \cdots \otimes  \fsquare(0.3cm,i_k) \otimes \cdots \otimes \fsquare(0.3cm,i_N)}

\font\germ=eufm10

\usepackage{graphicx,amsthm}
\usepackage{verbatim,amsmath,amscd,amssymb,color}
\begin{document}
\font\germ=eufm10
\def\bl{\bullet}
\def\aaa{@}
\title{\Large\bf RSK type correspondence of
Pictures and 
Littlewood-Richardson Crystals}
\vskip0.6cm
\author{Toshiki N\textsc{akashima}
\thanks{Department of Mathematics, 
Sophia University, Kioicho 7-1, Chiyoda-ku,
Tokyo 102-8554, Japan.\hfill\break
\qquad E-mail: toshiki{\aaa}mm.sophia.ac.jp\quad
:supported in part by  JSPS Grants in 
Aid for Scientific Research \#22540031.
}
\and Miki S\textsc{himojo}
\thanks{Department of Mathematics, 
Sophia University, Kioicho 7-1, Chiyoda-ku,
Tokyo 102-8554, Japan.\hfill\break
}}

\date{}
\maketitle


\abstract{We present a Robinson-Schensted-Knuth type 
one-to-one correspondence between 
the set of pictures and the set of pairs 
of  Littlewood-Richardson 
crystals.}

    


\newcommand{\cI}{{\mathcal I}}
\newcommand{\cA}{{\mathcal A}}
\newcommand{\cB}{{\mathcal B}}
\newcommand{\cC}{{\mathcal C}}
\newcommand{\cD}{{\mathcal D}}
\newcommand{\cF}{{\mathcal F}}
\newcommand{\cH}{{\mathcal H}}
\newcommand{\cK}{{\mathcal K}}
\newcommand{\cL}{{\mathcal L}}
\newcommand{\cM}{{\mathcal M}}
\newcommand{\cN}{{\mathcal N}}
\newcommand{\cO}{{\mathcal O}}
\newcommand{\cS}{{\mathcal S}}
\newcommand{\cV}{{\mathcal V}}
\newcommand{\fra}{\mathfrak a}
\newcommand{\frb}{\mathfrak b}
\newcommand{\frc}{\mathfrak c}
\newcommand{\frd}{\mathfrak d}
\newcommand{\fre}{\mathfrak e}
\newcommand{\frf}{\mathfrak f}
\newcommand{\frg}{\mathfrak g}
\newcommand{\frh}{\mathfrak h}
\newcommand{\fri}{\mathfrak i}
\newcommand{\frj}{\mathfrak j}
\newcommand{\frk}{\mathfrak k}
\newcommand{\frI}{\mathfrak I}
\newcommand{\fm}{\mathfrak m}
\newcommand{\frn}{\mathfrak n}
\newcommand{\frp}{\mathfrak p}
\newcommand{\fq}{\mathfrak q}
\newcommand{\frr}{\mathfrak r}
\newcommand{\frs}{\mathfrak s}
\newcommand{\frt}{\mathfrak t}
\newcommand{\fru}{\mathfrak u}
\newcommand{\frA}{\mathfrak A}
\newcommand{\frB}{\mathfrak B}
\newcommand{\frF}{\mathfrak F}
\newcommand{\frG}{\mathfrak G}
\newcommand{\frH}{\mathfrak H}
\newcommand{\frJ}{\mathfrak J}
\newcommand{\frN}{\mathfrak N}
\newcommand{\frP}{\mathfrak P}
\newcommand{\frT}{\mathfrak T}
\newcommand{\frU}{\mathfrak U}
\newcommand{\frV}{\mathfrak V}
\newcommand{\frX}{\mathfrak X}
\newcommand{\frY}{\mathfrak Y}
\newcommand{\frZ}{\mathfrak Z}
\newcommand{\rA}{\mathrm{A}}
\newcommand{\rC}{\mathrm{C}}
\newcommand{\rd}{\mathrm{d}}
\newcommand{\rB}{\mathrm{B}}
\newcommand{\rD}{\mathrm{D}}
\newcommand{\rE}{\mathrm{E}}
\newcommand{\rH}{\mathrm{H}}
\newcommand{\rK}{\mathrm{K}}
\newcommand{\rL}{\mathrm{L}}
\newcommand{\rM}{\mathrm{M}}
\newcommand{\rN}{\mathrm{N}}
\newcommand{\rR}{\mathrm{R}}
\newcommand{\rT}{\mathrm{T}}
\newcommand{\rZ}{\mathrm{Z}}
\newcommand{\bbA}{\mathbb A}
\newcommand{\bbC}{\mathbb C}
\newcommand{\bbG}{\mathbb G}
\newcommand{\bbF}{\mathbb F}
\newcommand{\bbH}{\mathbb H}
\newcommand{\bbP}{\mathbb P}
\newcommand{\bbN}{\mathbb N}
\newcommand{\bbQ}{\mathbb Q}
\newcommand{\bbR}{\mathbb R}
\newcommand{\bbV}{\mathbb V}
\newcommand{\bbZ}{\mathbb Z}
\newcommand{\adj}{\operatorname{adj}}
\newcommand{\Ad}{\mathrm{Ad}}
\newcommand{\Ann}{\mathrm{Ann}}
\newcommand{\rcris}{\mathrm{cris}}
\newcommand{\ch}{\mathrm{ch}}
\newcommand{\coker}{\mathrm{coker}}
\newcommand{\diag}{\mathrm{diag}}
\newcommand{\Diff}{\mathrm{Diff}}
\newcommand{\Dist}{\mathrm{Dist}}
\newcommand{\rDR}{\mathrm{DR}}
\newcommand{\ev}{\mathrm{ev}}
\newcommand{\Ext}{\mathrm{Ext}}
\newcommand{\cExt}{\mathcal{E}xt}
\newcommand{\fin}{\mathrm{fin}}
\newcommand{\Frac}{\mathrm{Frac}}
\newcommand{\GL}{\mathrm{GL}}
\newcommand{\Hom}{\mathrm{Hom}}
\newcommand{\hd}{\mathrm{hd}}
\newcommand{\rht}{\mathrm{ht}}
\newcommand{\id}{\mathrm{id}}
\newcommand{\im}{\mathrm{im}}
\newcommand{\inc}{\mathrm{inc}}
\newcommand{\ind}{\mathrm{ind}}
\newcommand{\coind}{\mathrm{coind}}
\newcommand{\Lie}{\mathrm{Lie}}
\newcommand{\Max}{\mathrm{Max}}
\newcommand{\mult}{\mathrm{mult}}
\newcommand{\op}{\mathrm{op}}
\newcommand{\ord}{\mathrm{ord}}
\newcommand{\pt}{\mathrm{pt}}
\newcommand{\qt}{\mathrm{qt}}
\newcommand{\rad}{\mathrm{rad}}
\newcommand{\res}{\mathrm{res}}
\newcommand{\rgt}{\mathrm{rgt}}
\newcommand{\rk}{\mathrm{rk}}
\newcommand{\SL}{\mathrm{SL}}
\newcommand{\soc}{\mathrm{soc}}
\newcommand{\Spec}{\mathrm{Spec}}
\newcommand{\St}{\mathrm{St}}
\newcommand{\supp}{\mathrm{supp}}
\newcommand{\Tor}{\mathrm{Tor}}
\newcommand{\Tr}{\mathrm{Tr}}
\newcommand{\wt}{\mathrm{wt}}
\newcommand{\Ab}{\mathbf{Ab}}
\newcommand{\Alg}{\mathbf{Alg}}
\newcommand{\Grp}{\mathbf{Grp}}
\newcommand{\Mod}{\mathbf{Mod}}
\newcommand{\Sch}{\mathbf{Sch}}\newcommand{\bfmod}{{\bf mod}}
\newcommand{\Qc}{\mathbf{Qc}}
\newcommand{\Rng}{\mathbf{Rng}}
\newcommand{\Top}{\mathbf{Top}}
\newcommand{\Var}{\mathbf{Var}}
\newcommand{\BB}{\mathbf{B}}
\newcommand{\gromega}{\langle\omega\rangle}
\newcommand{\lbr}{\begin{bmatrix}}
\newcommand{\rbr}{\end{bmatrix}}
\newcommand{\forb}{\bigcirc\kern-2.8ex \because}
\newcommand{\forbb}{\bigcirc\kern-3.0ex \because}
\newcommand{\forbbb}{\bigcirc\kern-3.1ex \because}
\newcommand{\cd}{commutative diagram }
\newcommand{\SpS}{spectral sequence}
\newcommand\C{\mathbb C}
\newcommand\hh{{\hat{H}}}
\newcommand\eh{{\hat{E}}}
\newcommand\F{\mathbb F}
\newcommand\fh{{\hat{F}}}

\def\AA{{\mathcal A}}
\def\al{\alpha}
\def\bq{B_q(\ge)}
\def\bqm{B_q^-(\ge)}
\def\bqz{B_q^0(\ge)}
\def\bqp{B_q^+(\ge)}
\def\beneme{\begin{enumerate}}
\def\beq{\begin{equation}}
\def\beqn{\begin{eqnarray}}
\def\beqnn{\begin{eqnarray*}}
\def\bigsl{{\hbox{\fontD \char'54}}}
\def\bbra#1,#2,#3{\left\{\begin{array}{c}\hspace{-5pt}
#1;#2\\ \hspace{-5pt}#3\end{array}\hspace{-5pt}\right\}}
\def\cd{\cdots}
\def\CC{\hbox{\bf C}}
\def\ddd{\hbox{\germ D}}
\def\del{\delta}
\def\Del{\Delta}
\def\Delr{\Delta^{(r)}}
\def\Dell{\Delta^{(l)}}
\def\Delb{\Delta^{(b)}}
\def\Deli{\Delta^{(i)}}
\def\Delre{\Delta^{\rm re}}
\def\douchik{\smash{\mathop{\sim }\limits^{k}}}
\def\douchic{\smash{\mathop{\sim}\limits^{c}}}

\def\ei{e_i}
\def\eit{\tilde{e}_i}
\def\eneme{\end{enumerate}}
\def\ep{\epsilon}
\def\eeq{\end{equation}}
\def\eeqn{\end{eqnarray}}
\def\eeqnn{\end{eqnarray*}}
\def\fit{\tilde{f}_i}
\def\FF{{\rm F}}
\def\ft{\tilde{f}}
\def\gau#1,#2{\left[\begin{array}{c}\hspace{-5pt}#1\\
\hspace{-5pt}#2\end{array}\hspace{-5pt}\right]}
\def\ge{\hbox{\germ g}}
\def\gl{\hbox{\germ gl}}
\def\hom{{\hbox{Hom}}}
\def\ify{\infty}
\def\io{\iota}
\def\kp{k^{(+)}}
\def\km{k^{(-)}}
\def\llra{\relbar\joinrel\relbar\joinrel\relbar\joinrel\rightarrow}
\def\lan{\langle}
\def\lar{\longrightarrow}
\def\max{{\rm max}}
\def\ME{\rm ME}
\def\lm{\lambda}
\def\Lm{\Lambda}
\def\mapright#1{\smash{\mathop{\longrightarrow}\limits^{#1}}}
\def\mm{{\bf{\rm m}}}
\def\nd{\noindent}
\def\nn{\nonumber}
\def\nnn{\hbox{\germ n}}
\def\catob{{\mathcal O}(B)}
\def\oint{{\mathcal O}_{\rm int}(\ge)}
\def\ot{\otimes}
\def\op{\oplus}
\def\opi{\ovl\pi_{\lm}}
\def\ovl{\overline}
\def\plm{\Psi^{(\lm)}_{\io}}
\def\qq{\qquad}
\def\q{\quad}
\def\qed{\hfill\framebox[2mm]{}}
\def\QQ{\hbox{\bf Q}}
\def\qi{q_i}
\def\qii{q_i^{-1}}
\def\ra{\rightarrow}
\def\ran{\rangle}
\def\rlm{r_{\lm}}
\def\ssl{\mathfrak{sl}}
\def\slh{\widehat{\ssl_2}}
\def\ge{\hbox{\germ g}}
\def\ti{t_i}
\def\tii{t_i^{-1}}
\def\til{\tilde}
\def\tm{\times}
\def\tt{{\hbox{\germ{t}}}}
\def\ttt{\hbox{\germ t}}
\def\ua{U_{\AA}}
\def\ue{U_{\vep}}
\def\uq{U_q(\ge)}
\def\ufin{U^{\rm fin}_{\vep}}
\def\ufinp{(U^{\rm fin}_{\vep})^+}
\def\ufinm{(U^{\rm fin}_{\vep})^-}
\def\ufinz{(U^{\rm fin}_{\vep})^0}
\def\uqm{U^-_q(\ge)}
\def\uqp{U^+_q(\ge)}
\def\uqmq{{U^-_q(\ge)}_{\bf Q}}
\def\uqpm{U^{\pm}_q(\ge)}
\def\uqq{U_{\bf Q}^-(\ge)}
\def\uqz{U^-_{\bf Z}(\ge)}
\def\ures{U^{\rm res}_{\AA}}
\def\urese{U^{\rm res}_{\vep}}
\def\uresez{U^{\rm res}_{\vep,\ZZ}}
\def\util{\widetilde\uq}
\def\uup{U^{\geq}}
\def\ulow{U^{\leq}}
\def\bup{B^{\geq}}
\def\blow{\ovl B^{\leq}}
\def\vep{\varepsilon}
\def\vp{\varphi}
\def\vpi{\varphi^{-1}}
\def\VV{{\mathcal V}}
\def\xii{\xi^{(i)}}
\def\Xiioi{\Xi_{\io}^{(i)}}
\def\WW{{\mathcal W}}
\def\wtil{\widetilde}
\def\what{\widehat}
\def\wpi{\widehat\pi_{\lm}}
\def\ZZ{\mathbb Z}
\def\spsp(#1,#2){\begin{pmatrix}
#1,\\ \hline#2\end{pmatrix}}

\theoremstyle{definition}
\newtheorem{df}{\bf Definition}[section]
\newtheorem{pro}[df]{\bf Proposition}
\newtheorem{thm}[df]{\bf Theorem}
\newtheorem{lem}[df]{\bf Lemma}
\newtheorem{ex}[df]{\bf Example}
\newtheorem{cor}[df]{\bf Corollary}
\newtheorem{con}[df]{Conjecture}

\def\m@th{\mathsurround=0pt}

\def\fsquare(#1,#2){
\hbox{\vrule$\hskip-0.4pt\vcenter to #1{\normalbaselines\m@th
\hrule\vfil\hbox to #1{\hfill$\scriptstyle #2$\hfill}\vfil\hrule}$\hskip-0.4pt
\vrule}}

\def\addsquare(#1,#2){\hbox{$
	\dimen1=#1 \advance\dimen1 by -0.8pt
	\vcenter to #1{\hrule height0.4pt depth0.0pt%
	\hbox to #1{%
	\vbox to \dimen1{\vss%
	\hbox to \dimen1{\hss$\scriptstyle~#2~$\hss}%
	\vss}%
	\vrule width0.4pt}%
	\hrule height0.4pt depth0.0pt}$}}

\def\Fsquare(#1,#2){
\hbox{\vrule$\hskip-0.4pt\vcenter to #1{\normalbaselines\m@th
\hrule\vfil\hbox to #1{\hfill$#2$\hfill}\vfil\hrule}$\hskip-0.4pt
\vrule}}

\def\Addsquare(#1,#2){\hbox{$
	\dimen1=#1 \advance\dimen1 by -0.8pt
	\vcenter to #1{\hrule height0.4pt depth0.0pt%
	\hbox to #1{%
	\vbox to \dimen1{\vss%
	\hbox to \dimen1{\hss$~#2~$\hss}%
	\vss}%
	\vrule width0.4pt}%
	\hrule height0.4pt depth0.0pt}$}}

\def\hfourbox(#1,#2,#3,#4){%
	\fsquare(0.3cm,#1)\addsquare(0.3cm,#2)\addsquare(0.3cm,#3)\addsquare(0.3cm,#4)}

\def\Hfourbox(#1,#2,#3,#4){%
	\Fsquare(0.4cm,#1)\Addsquare(0.4cm,#2)\Addsquare(0.4cm,#3)\Addsquare(0.4cm,#4)}

\def\HHfourbox(#1,#2,#3,#4){%
	\Fsquare(0.8cm,#1)\Addsquare(0.8cm,#2)\Addsquare(0.8cm,#3)\Addsquare(0.8cm,#4)}

\def\fsq(#1){%
          \fsquare(0.3cm,#1)}

\def\hthreebox(#1,#2,#3){%
	\fsquare(0.3cm,#1)\addsquare(0.3cm,#2)\addsquare(0.3cm,#3)}

\def\htwobox(#1,#2){%
	\fsquare(0.3cm,#1)\addsquare(0.3cm,#2)}

\def\vfourbox(#1,#2,#3,#4){%
	\normalbaselines\m@th\offinterlineskip
	\vcenter{\hbox{\fsquare(0.3cm,#1)}
	      \vskip-0.4pt
	      \hbox{\fsquare(0.3cm,#2)}	
	      \vskip-0.4pt
	      \hbox{\fsquare(0.3cm,#3)}	
	      \vskip-0.4pt
	      \hbox{\fsquare(0.3cm,#4)}}}

\def\VVfourbox(#1,#2,#3,#4){%
	\normalbaselines\m@th\offinterlineskip
	\vcenter{\hbox{\Fsquare(0.8cm,#1)}
	      \vskip-0.4pt
	      \hbox{\Fsquare(0.8cm,#2)}	
	      \vskip-0.4pt
	      \hbox{\Fsquare(0.8cm,#3)}	
	      \vskip-0.4pt
	      \hbox{\Fsquare(0.8cm,#4)}}}

\def\Vfourbox(#1,#2,#3,#4){%
	\normalbaselines\m@th\offinterlineskip
	\vcenter{\hbox{\Fsquare(0.4cm,#1)}
	      \vskip-0.4pt
	      \hbox{\Fsquare(0.4cm,#2)}	
	      \vskip-0.4pt
	      \hbox{\Fsquare(0.4cm,#3)}	
	      \vskip-0.4pt
	      \hbox{\Fsquare(0.4cm,#4)}}}

\def\vthreebox(#1,#2,#3){%
	\normalbaselines\m@th\offinterlineskip
	\vcenter{\hbox{\fsquare(0.3cm,#1)}
	      \vskip-0.4pt
	      \hbox{\fsquare(0.3cm,#2)}	
	      \vskip-0.4pt
	      \hbox{\fsquare(0.3cm,#3)}}}

\def\vtwobox(#1,#2){%
	\normalbaselines\m@th\offinterlineskip
	\vcenter{\sbox{\fsquare(0.3cm,#1)}
	      \vskip-0.4pt
	      \hbox{\fsquare(0.3cm,#2)}}}

\def\vtwobox2(#1,#2){%
	\normalbaselines\m@th\offinterlineskip
	\vcenter{\hbox{#1}
	      \vskip-0.4pt
	      \hbox{\fsquare(0.3cm,#2)}}}

\def\Hthreebox(#1,#2,#3){%
	\Fsquare(0.4cm,#1)\Addsquare(0.4cm,#2)\Addsquare(0.4cm,#3)}

\def\HHthreebox(#1,#2,#3){%
	\Fsquare(0.8cm,#1)\Addsquare(0.8cm,#2)\Addsquare(0.8cm,#3)}

\def\Htwobox(#1,#2){%
	\Fsquare(0.4cm,#1)\Addsquare(0.4cm,#2)}

\def\H6twobox(#1,#2){%
	\Fsquare(0.6cm,#1)\Addsquare(0.6cm,#2)}

\def\HHtwobox(#1,#2){%
	\Fsquare(0.8cm,#1)\Addsquare(0.8cm,#2)}

\def\Vthreebox(#1,#2,#3){%
	\normalbaselines\m@th\offinterlineskip
	\vcenter{\hbox{\Fsquare(0.4cm,#1)}
	      \vskip-0.4pt
	      \hbox{\Fsquare(0.4cm,#2)}	
	      \vskip-0.4pt
	      \hbox{\Fsquare(0.4cm,#3)}}}

\def\Vtwobox(#1,#2){%
	\normalbaselines\m@th\offinterlineskip
	\vcenter{\hbox{\Fsquare(0.4cm,#1)}
	      \vskip-0.4pt
	      \hbox{\Fsquare(0.4cm,#2)}}}

\def\twoone(#1,#2,#3){%
	\normalbaselines\m@th\offinterlineskip
	\vcenter{\hbox{\htwobox({#1},{#2})}
	      \vskip-0.4pt
	      \hbox{\fsquare(0.3cm,#3)}}}

\def\twothree(#1,#2,#3,#4,#5){%
	\normalbaselines\m@th\offinterlineskip
	\vcenter{\hbox{\htwobox({#1},{#2})}
	      \vskip-0.4pt
	      \hbox{\hthreebox({#3},{#4},{#5})}}}

\def\threethree(#1,#2,#3,#4,#5,#6){%
	\normalbaselines\m@th\offinterlineskip
	\vcenter{\hbox{\hthreebox({#1},{#2},{#3})}
	      \vskip-0.4pt
	      \hbox{\hthreebox({#4},{#5},{#6})}}}

\def\threethreeone(#1,#2,#3,#4,#5,#6,#7){%
	\normalbaselines\m@th\offinterlineskip
	\vcenter{\hbox{\hthreebox({#1},{#2},{#3})}
	      \vskip-0.4pt
	      \hbox{\hthreebox({#4},{#5},{#6})}
              \vskip-0.4pt
              \hbox{\fsquare(0.3cm,#7)}}}

\def\threethreetwo(#1,#2,#3,#4,#5,#6,#7,#8){%
	\normalbaselines\m@th\offinterlineskip
	\vcenter{\hbox{\hthreebox({#1},{#2},{#3})}
	      \vskip-0.4pt
	      \hbox{\hthreebox({#4},{#5},{#6})}
              \vskip-0.4pt
              \hbox{\htwobox({#7},{#8})}}}

\def\Twoone(#1,#2,#3){%
	\normalbaselines\m@th\offinterlineskip
	\vcenter{\hbox{\Htwobox({#1},{#2})}
	      \vskip-0.4pt
	      \hbox{\Fsquare(0.4cm,#3)}}}

\def\TTwoone(#1,#2,#3){%
	\normalbaselines\m@th\offinterlineskip
	\vcenter{\hbox{\H6twobox({#1},{#2})}
	      \vskip-0.4pt
	      \hbox{\Fsquare(0.6cm,#3)}}}

\def\threeone(#1,#2,#3,#4){%
	\normalbaselines\m@th\offinterlineskip
	\vcenter{\hbox{\hthreebox({#1},{#2},{#3})}
	      \vskip-0.4pt
	      \hbox{\fsquare(0.3cm,#4)}}}

\def\Threeone(#1,#2,#3,#4){%
	\normalbaselines\m@th\offinterlineskip
	\vcenter{\hbox{\Hthreebox({#1},{#2},{#3})}
	      \vskip-0.4pt
	      \hbox{\Fsquare(0.4cm,#4)}}}

\def\Threetwo(#1,#2,#3,#4,#5){%
	\normalbaselines\m@th\offinterlineskip
	\vcenter{\hbox{\Hthreebox({#1},{#2},{#3})}
	      \vskip-0.4pt
	      \hbox{\Htwobox({#4},{#5})}}}

\def\threetwo(#1,#2,#3,#4,#5){%
	\normalbaselines\m@th\offinterlineskip
	\vcenter{\hbox{\hthreebox({#1},{#2},{#3})}
	      \vskip-0.4pt
	      \hbox{\htwobox({#4},{#5})}}}

\def\twotwo(#1,#2,#3,#4){%
	\normalbaselines\m@th\offinterlineskip
	\vcenter{\hbox{\htwobox({#1},{#2})}
	      \vskip-0.4pt
	      \hbox{\htwobox({#3},{#4})}}}

\def\Twotwo(#1,#2,#3,#4){%
	\normalbaselines\m@th\offinterlineskip
	\vcenter{\hbox{\Htwobox({#1},{#2})}
	      \vskip-0.4pt
	      \hbox{\Htwobox({#3},{#4})}}}

\def\TTwotwo(#1,#2,#3,#4){%
	\normalbaselines\m@th\offinterlineskip
	\vcenter{\hbox{\H6twobox({#1},{#2})}
	      \vskip-0.4pt
	      \hbox{\H6twobox({#3},{#4})}}}

\def\twooneone(#1,#2,#3,#4){%
	\normalbaselines\m@th\offinterlineskip
	\vcenter{\hbox{\htwobox({#1},{#2})}
	      \vskip-0.4pt
	      \hbox{\fsquare(0.3cm,#3)}
	      \vskip-0.4pt
	      \hbox{\fsquare(0.3cm,#4)}}}

\def\Twooneone(#1,#2,#3,#4){%
	\normalbaselines\m@th\offinterlineskip
	\vcenter{\hbox{\Htwobox({#1},{#2})}
	      \vskip-0.4pt
	      \hbox{\Fsquare(0.4cm,#3)}
	      \vskip-0.4pt
	      \hbox{\Fsquare(0.4cm,#4)}}}

\def\Twotwoone(#1,#2,#3,#4,#5){%
	\normalbaselines\m@th\offinterlineskip
	\vcenter{\hbox{\Htwobox({#1},{#2})}
	      \vskip-0.4pt
	      \hbox{\Htwobox({#3},{#4})}
              \vskip-0.4pt
	      \hbox{\Fsquare(0.4cm,#5)}}}

\def\twotwoone(#1,#2,#3,#4,#5){%
	\normalbaselines\m@th\offinterlineskip
	\vcenter{\hbox{\htwobox({#1},{#2})}
	      \vskip-0.4pt
	      \hbox{\htwobox({#3},{#4})}
              \vskip-0.4pt
	      \hbox{\fsquare(0.3cm,#5)}}}

\def\twotwotwo(#1,#2,#3,#4,#5,#6){%
	\normalbaselines\m@th\offinterlineskip
	\vcenter{\hbox{\htwobox({#1},{#2})}
	      \vskip-0.4pt
	      \hbox{\htwobox({#3},{#4})}
              \vskip-0.4pt
	      \hbox{\htwobox({#5},{#6})}}}

\def\threetwoone(#1,#2,#3,#4,#5,#6){%
	\normalbaselines\m@th\offinterlineskip
	\vcenter{\hbox{\hthreebox({#1},{#2},{#3})}
	      \vskip-0.4pt
	      \hbox{\htwobox({#4},{#5})}
              \vskip-0.4pt
	      \hbox{\fsquare(0.3cm,#6)}}}

\def\fourthreetwoone{%
	\normalbaselines\m@th\offinterlineskip
	\vcenter{\hbox{\hfourbox(,,,)}
	      \vskip-0.4pt
              \hbox{\hthreebox(,,)}
              \vskip-0.4pt
	      \hbox{\htwobox(,)}
              \vskip-0.4pt
	      \hbox{\fsquare(0.3cm,)}}}

\def\fourtwoone(#1,#2,#3,#4,#5,#6,#7){%
	\normalbaselines\m@th\offinterlineskip
	\vcenter{\hbox{\hfourbox({#1},{#2},{#3},{#4})}
	      \vskip-0.4pt
	      \hbox{\htwobox({#5},{#6})}
              \vskip-0.4pt
	      \hbox{\fsquare(0.3cm,{#7})}}}

\def\fourthreetwo(#1,#2,#3,#4,#5,#6,#7,#8,#9){%
	\normalbaselines\m@th\offinterlineskip
	\vcenter{\hbox{\hfourbox({#1},{#2},{#3},{#4})}
	      \vskip-0.4pt
              \hbox{\hthreebox({#5},{#6},{#7})}
              \vskip-0.4pt
	      \hbox{\htwobox({#8},{#9})}}}
\def\fourthreeone(#1,#2,#3,#4,#5,#6,#7,#8){%
	\normalbaselines\m@th\offinterlineskip
	\vcenter{\hbox{\hfourbox({#1},{#2},{#3},{#4})}
	      \vskip-0.4pt
              \hbox{\hthreebox({#5},{#6},{#7})}
              \vskip-0.4pt
	      \hbox{\fsquare(0.3cm,{#8})}}}

\def\TThreetwoone(#1,#2,#3,#4,#5,#6){%
	\normalbaselines\m@th\offinterlineskip
	\vcenter{\hbox{\HHthreebox({#1},{#2},{#3})}
	      \vskip-0.4pt
	      \hbox{\HHtwobox({#4},{#5})}
              \vskip-0.4pt
	      \hbox{\fsquare(0.8cm,#6)}}}

\def\Threeoneone(#1,#2,#3,#4,#5){%
	\normalbaselines\m@th\offinterlineskip
	\vcenter{\hbox{\Hthreebox({#1},{#2},{#3})}
	      \vskip-0.4pt
	      \hbox{\Fsquare(0.4cm,#4)}
              \vskip-0.4pt
	      \hbox{\Fsquare(0.4cm,#5)}}}

\def\threeoneone(#1,#2,#3,#4,#5){%
	\normalbaselines\m@th\offinterlineskip
	\vcenter{\hbox{\hthreebox({#1},{#2},{#3})}
	      \vskip-0.4pt
	      \hbox{\fsquare(0.3cm,#4)}
              \vskip-0.4pt
	      \hbox{\fsquare(0.3cm,#5)}}}

\def\FFourtwoone(#1,#2,#3,#4,#5,#6,#7){%
	\normalbaselines\m@th\offinterlineskip
	\vcenter{\hbox{\HHfourbox({#1},{#2},{#3},{#4})}
	      \vskip-0.4pt
	      \hbox{\HHtwobox({#5},{#6})}
              \vskip-0.4pt
	      \hbox{\fsquare(0.8cm,#7)}}}

\def\a{\fsquare(0.3cm){1}\addsquare(0.3cm)(2)\addsquare(0.3cm)(3)}

\def\b{\hbox{%
	\normalbaselines\m@th\offinterlineskip
	\vcenter{\hbox{\fsquare(0.3cm){2}}\vskip-0.4pt\hbox{\fsquare(0.3cm){2}}}}}

\def\c{\hbox{\normalbaselines\m@th\offinterlineskip%
	\vcenter{\hbox{\a}\vskip-0.4pt\hbox{\b}}}}


\dimen1=0.4cm\advance\dimen1 by -0.8pt
\def\ffsquare#1{%
	\fsquare(0.4cm,\hbox{#1})}

\def\naga{%
	\hbox{$\vcenter to 0.4cm{\normalbaselines\m@th
	\hrule\vfil\hbox to 1.2cm{\hfill$\cdots$\hfill}\vfil\hrule}$}}

\def\vnaga{\normalbaselines\m@th\baselineskip0pt\offinterlineskip%
	\vrule\vbox to 1.2cm{\vskip7pt\hbox to \dimen1{$\hfil\vdots\hfil$}\vfil}\vrule}

\def\dhbox{\hbox{$\ffsquare 1 \naga \ffsquare N$}}

\def\dvbox{\hbox{\normalbaselines\m@th\baselineskip0pt\offinterlineskip\vbox{%
	  \hbox{$\ffsquare 1$}\vskip-0.4pt\hbox{$\vnaga$}\vskip-0.4pt\hbox{$\ffsquare N$}}}}

\def\sq(#1){\fsquare(0.4cm,#1)}
\def\Sq(#1){\fsquare(0.5cm,#1)}
\def\SSq(#1){\fsquare(0.9cm,#1)}
\def\Sqj{\Sq(j)}
\def\Sqjb{\Sq(\ovl j)}
\def\Sqjm{\Sq(\scriptstyle{j-1})}
\def\Sqjp{\Sq(\scriptstyle{j+1})}
\def\Sqjmb{\Sq(\scriptstyle{\ovl{j-1}})}
\def\SSqj{\SSq(j)}
\def\Sqjpb{\Sq(\scriptstyle{\ovl{j+1}})}
\def\SSqjb{\SSq(\ovl j)}
\def\SSqjm{\SSq(\scriptstyle{j-1})}
\def\SSqjp{\SSq(\scriptstyle{j+1})}
\def\SSqjmb{\SSq(\scriptstyle{\ovl{j-1}})}
\def\SSqjpb{\SSq(\scriptstyle{\ovl{j+1}})}

\def\mapright#1{\smash{\mathop{\longrightarrow}\limits^{#1}}}
\def\map#1{\smash{\mathop{\longmapsto}\limits^{#1}}}
\def\maplr#1{\smash{\mathop{\longleftrightarrow}\limits^{#1}}}
\def\Longsquare(#1,#2,#3){
\hbox{\vrule$\hskip-0.4pt\vcenter to #1{\normalbaselines\m@th
\hrule\vfil\hbox to #2{\hfill$#3$\hfill}\vfil\hrule}$\hskip-0.4pt
\vrule}}

\renewcommand{\thesection}{\arabic{section}}
\section{Introduction}
\setcounter{equation}{0}
\renewcommand{\theequation}{\thesection.\arabic{equation}}

Combinatorics of pictures has been initiated in 
\cite{CS,CS2,JP,Z}. 
Picture is a certain bijective order
morphism
between two skew Young diagrams with some partial/total orders.
The remarkable result for pictures is that 
there exists a kind of RSK type one to one correspondence
as follows.
Let $\kappa^i$ $(i=1,2)$ be skew Young diagrams with 
$|\kappa^1|=|\kappa^2|(=N)$. There exists a bijection:
\begin{equation}
 {\mathbf P}(\kappa^1,\kappa^2)\q\maplr{1:1}\q
 \coprod_{\mu}
 \left(
{\mathbf P}(\mu,\kappa^1)
 \times  {\mathbf P}(\mu,\kappa^2)\right),
\label{ori}
\end{equation}
where $\mu$ runs over the set of Young diagrams with $|\mu|=N$
and $ {\mathbf P}(\kappa^1,\kappa^2)$ is a set of pictures
from $\kappa^1$ to $\kappa^2$.
Since some set of pictures can be identified with a set of 
permutations, this correspondence can be seen as
an analogue of the RSK correspondence.
In \cite{FG,vL}, certain generalizations have been done using 
various combinatorial methods.

In \cite{NS,NS2}, we introduced the one to one correspondence
between ''Littlewood-Richardson crystals'' and  pictures.
\begin{equation}
 {\mathbf P}(\mu,\nu\setminus\lm)\q\maplr{1:1}\q
\mathbf{B}(\mu)^{\nu}_{\lambda},\label{1iso}
\end{equation}
where $\lm,\mu,\nu$ are Young diagrams with $|\lm|+|\mu|=|\nu|$.
This seems to give a new interpretation of pictures from the
view point of theory of crystal bases.

In this article, we shall describe the following bijections
\begin{equation}
{\mathbf P}(\kappa^1,\kappa^2)\q\maplr{1:1}\q
{\mathbf S}(\kappa^1,\kappa^2)\q\maplr{1:1}\q
{\mathbf W}(\kappa^1,\kappa^2)\q\maplr{1:1}\q
 \coprod_{\mu} \left(\mathbf{B}(\mu)^{\nu^1}_{\lambda^1}
 \times \mathbf{B}(\mu)^{\nu^2}_{\lambda^2}\right),
\label{4iso}
\end{equation}
where ${\mathbf P}(\kappa^1,\kappa^2)$ is a set of pictures 
from $\kappa^1$ to $\kappa^2$, 
${\mathbf S}(\kappa^1,\kappa^2)$ is a set of 
Littlewood-Richardson skew tableaux associated with 
$(\kappa^1,\kappa^2)$, 
${\mathbf W}(\kappa^1,\kappa^2)$ is a set of lexicographic 
two-rowed array (of column type)
associated with $(\kappa^1,\kappa^2)$ and 
the last one is a set of pairs of Littlewood-Richardson crystals.
Thus, applying (\ref{1iso}) to the last one in (\ref{4iso})
we obtain the original correspondence (\ref{ori}).
The pictures treated in this article are defined by the 
order $J$ (see Sect.2), 
which is a kind of admissible orders. More 
general setting, namely defined by general admissible orders
will be discussed elsewhere.

As is well known that the crystal ${\mathbf B}(\mu)$ of type
${\mathfrak gl}_n$ (or ${\mathfrak sl}_n$) is realized
as the set of Young tableaux \cite{KN} and 
the Littlewood-Richardson crystal ${\mathbf B}(\mu)^\nu_\lm$ 
is a subset of ${\mathbf B}(\mu)$ with the
certain special conditions 'highest conditions' 
\cite{N1,NS,NS2}. Thus, the last term in (\ref{4iso})
is a set of pairs of same shaped Young tableaux and then 
bijections in (\ref{4iso}) turn out 
to be a generalization of the
RSK correspondence.

As claimed in \cite{NS,NS2}, these methods would 
open the door to generalize the theory of pictures to 
wider classes. Indeed, 
in preparing this manuscript, we received the preprint
'Admissible pictures and $U_q(gl(m; n))$-
Littlewood-Richardson tableaux' by J.I.Hung, 
S-J.Kang and Y-W. Lyoo, which gives the first bijection 
in (\ref{4iso}) and generalizes it to the 
the super case $U_q(gl(m; n))$. This is a kind of the
evidence of our claims, unfortunately, which was not done by us.

The organizations of the article is as follows:
in Sect.2 and 3, the basics of pictures and crystals are 
reviewed.  In Sect.4, we introduce several combinatorial
procedures and notions required in this article; 
column bumping, RSK correspondence,
Knuth equivalence, crystal equivalence and etc.
The main theorem is given in Sect.5. and its proof is described 
separately in the subsequent sections.

\renewcommand{\thesection}{\arabic{section}}
\section{Pictures}
\setcounter{equation}{0}
\renewcommand{\theequation}{\thesection.\arabic{equation}}

\subsection{Young diagrams and Young tableaux}\label{pT}
Let $\lambda=(\lambda_1,\lambda_2,\cdots,\lambda_m)$
be a Young diagram or a partition, which satisfies
$\lambda_1 \geq\lambda_2 \geq \cdots \geq \lambda_m \geq 0$.
Let $\lambda$ and $\mu$ be Young diagrams with $\mu\subset\lm$.
A {\it skew diagram} $\lambda \setminus \mu$ is obtained by 
subtracting set-theoretically $\mu$ from $\lm$.

In this article we frequently consider 
a (skew) Young diagram as a subset of $\bbN\times\bbN$ by 
identifying the box in the $i$-th row and the $j$-th column 
with $(i,j)\in \bbN\times\bbN$.

\begin{ex}@A Young diagram $\lambda=(2,2,1)$ is expressed by 
$\left\{
(1,1), (1,2),
(2,1), (2,2),
(3,1)
\right\}$.
\end{ex}

As in \cite{F}, 
in the sequel,
a ''(skew) Young tableau'' means  a semi-standard (skew) tableau.
For a skew Young tableau $S$ of shape $\lm\setminus \mu$, 
we also consider a ''coordinate '' in $\bbN\times\bbN$
like as a skew diagram $\lm\setminus\mu$. 
Then an entry of $S$ in $(i,j)$ is denoted by $S_{i,j}$ and 
called $(i,j)$-entry.
For $k>0$, define (\cite{NS})
\begin{equation}
\label{Y-diag}
S^{(k)} = \{ (l,m) \in \lm\setminus\mu | S_{l,m} = k \}.
\end{equation}
There is no two elements in one column in $S^{(k)}$. 
For a skew Young tableau $S$ with $(i,j)$-entry $S_{i,j}=k$, 
we define $p(S;i,j)$  (\cite{NS}) as
the number of $(i,j)$-entry 
from the right in $S^{(k)}$.

\subsection{Picture}\label{pic}

First, we shall introduce the original notion of 
''picture'' as in \cite{Z}.

We define the following two kinds of orders on
a subset $X\subset\bbN \times \bbN$: For
$(a,b), \,\,(c,d)\in X$,
\begin{enumerate}
\item $(a,b)\leqslant_P(c,d)$ iff
$a\leq c \text{ and } b \leq d$.
\item $(a,b) \leqslant_J (c,d)$ iff
$a < c \text{, or }  a = c \text{ and } b\geq d$.
\end{enumerate}
Note that the order $\leqslant_P$ is a partial order 
and $\leqslant_J$ is a total order.
\begin{df}[\cite{Z}] 
Let $X,Y \subset \bbN \times \bbN$.
\begin{enumerate}
\item
A map $f : X \to Y$ is said to be 
{\it PJ-standard} if it satisfies
\[
\text{For }(a,b),(c,d) \in X, 
\text{if }(a,b) \leqslant_P (c,d), \text{ then }  
f(a,b) \leqslant_J f(c,d).
\]
\item
A map $f:X\to Y$ is a {\it picture} if it is 
bijective and
both $f$ and $f^{-1}$ are PJ-standard.
\end{enumerate}
\end{df}

Taking two skew Young diagrams
$\kappa^1,\kappa^2\subset \bbN\times\bbN$, 
denote the set of pictures by:
\[
{\mathbf P}(\kappa^1,\kappa^2):= \{ f : \kappa^1 \to 
\kappa^2 \,|\,f\text{ is a picture.}\}
\]

Next, we shall generalize the notion of pictures by using
a total order on a subset 
$X\subset \bbN\times \bbN$, called an
``{\it admissible order}'', though we do not treat this 
generalization in this article:
\begin{df}
\begin{enumerate}
\item
A total order $\leqslant_A$ on $X\subset\bbN\times\bbN$ 
is called {\it admissible}
if it satisfies:
\[
\text{For any }(a,b),\,\,(c,d)\in X\text{ if }
a\leq c\text{ and }b\geq d\text{ then } (a,b)\leqslant_A (c,d).
\]
\item
For $X,Y\subset\bbN\times\bbN$ and a map $f:X\to Y$, if $f$ satisfies
that if $(a,b)\leqslant_P(c,d)$, then $f(a,b)\leqslant_A f(c,d)$ for 
any $(a,b),~(c,d)\in X$, then $f$ is called $PA$-standard.
\item
Let $\leqslant_A$ (resp. $\leqslant_{A'}$) be an admissible order on 
$X\text{(resp. }Y)\subset\bbN\times\bbN$. A bijective map 
$f:X\to Y$ is called an $(A,A')$-{\it admissible picture} or 
simply, an {\it admissible picture} if $f$ is $PA$-standard and 
$f^{-1}$ is $PA'$-standard.
\end{enumerate}
\end{df}

\renewcommand{\thesection}{\arabic{section}}
\section{Crystals}
\setcounter{equation}{0}
\renewcommand{\theequation}{\thesection.\arabic{equation}}

The basic references for the theory of crystals are
\cite{K1},\cite{K2}.
\subsection{Readings and Additions}
Let $\BB=\{\fsq(i)\,|\,1\leq i\leq n+1\}$ 
be the crystal of the vector representation 
$V(\Lm_1)$ of the quantum group $U_q(A_n)$ (\cite{KN}).
As in \cite{NS}, we shall identify a dominant integral weight 
of type $A_n$ with a Young diagram in the standard way, 
{\it e.g.,} the fundamental weight $\Lm_1$ is identified with 
a square box $\fsq()$. For a Young diagram $\lm$, 
let $B(\lm)$ be the crystal of the finite-dimensional 
irreducible $U_q(A_n)$-module $V(\lm)$.
Set $N:=|\lm|$. Then there exists an embedding of crystals:
$B(\lm)\hookrightarrow \BB^N$ and an element in $B(\lm)$ 
is realized by a Young tableau of shape $\lm$
(\cite{KN}). 
Note that these embedding can be extended to skew tableaux, that is, 
there exists a embedding of crystals
$S(\kappa)\hookrightarrow {\mathbf B}^{\ot N}$, 
where $S(\kappa)$ is the set of skew tableaux of shape $\kappa$
$N=|\kappa|$ (\cite{HK}). Indeed, 
$S(\kappa)$ is a direct sum of certain $B(\lm)$'s.
Such an embedding is not unique,
which is called a 'reading' and described by:
\begin{df}[\cite{HK}]
Let $A$ be an admissible order on 
a (skew) Young diagram $\lm$  with $|\lm|=N$.
For $T\in B(\lm)$ (resp. $S(\lm)$), by reading the entries 
in $T$ according to $A$, 
we obtain the map 
\[
 R_A:B(\lm) (\text{resp. }S(\lm))
\longrightarrow {\mathbf B}^{\ot N}\q (T\mapsto
\fsq(i_1)\ot\cd\ot\fsq(i_N))),
\]
which is called an {\it admissible reading}
associated with the order $A$.
The map $R_A$ is 
an embedding of crystals. In particular, in case that taking the
order $J$ as an admissible order, we denote the embedding $R_J$
by ME and call it a {\it middle-eastern reading}.
\end{df}

\begin{df}
For $i\in\{1,2,\cd,n+1\}$ and a Young diagram 
$\lambda=(\lambda_1,\lambda_2,\cdots, \lambda_n)$, we define
\[
\lambda[i]:=
(\lambda_1,\lambda_2,\cdots,\lambda_i+1,\cdots, \lambda_n)
\]
which is said to be an {\it addition} of $i$ to $\lm$.
In general, 
for $i_1,i_2,\cd,i_N\in\{1,2,\cd,n+1\} $ and a Young diagram 
$\lm$, we define 
\[
\lambda[i_1,i_2,\cdots,i_N]
:=(\cd ((\lambda[i_1])[i_2])\cd)[i_N],
\]
which is called an {\it addition} of $i_1,\cd,i_N$ to $\lm$.
\end{df}
\begin{ex}
For a sequence ${\bf i}=31212$, 
the addition of $\bf i$ to $\lambda = \twoone(,,)$ is:
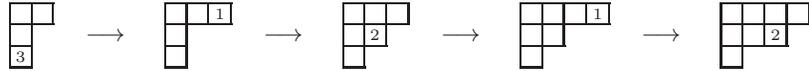
\[
\twooneone(,,,3) \quad \longrightarrow \quad \threeoneone(,,1,,)
\quad \longrightarrow \quad \threetwoone(,,,,2,)
\quad \longrightarrow \quad \fourtwoone(,,,1,,,) \quad \longrightarrow
 \quad \fourthreeone(,,,,,,2,).
\]
\end{ex}
\vspace{0.5cm}
{\sl Remark.}
For a Young diagram $\lm$, 
an addition $\lm[i_1,\cd,i_N]$ is not 
necessarily a Young diagram.
For instance, a sequence ${\bf i'}=22133$ and 
$\lambda=(2,2)$,
the addition $\lm[{\bf i'}]=(3,3,2)$ is a Young diagram.
But, in the second step of the addition, it becomes the diagram 
$\lm[2,2]=(2,3)$, which is not a Young diagram.

\subsection{Littlewood-Richardson Crystal}\label{subsec:LR}
As an application of the description of 
crystal bases of type $A_n$,
we see so-called ``Littlewood-Richardson rule'' of type $A_n$.

For a sequence ${\mathbf i}=
i_1,i_2,\cd,i_N\,\,(i_j\in \{1,2,\cd,n+1\})$ and 
a Young diagram $\lm$, let $\til\lm:=\lm[i_1,i_2,\cd,i_{N}]$ be
an addition of $i_1,i_2,\cd,i_N$ to $\lm$. Then set
\[
{\mathbf B}(\lm:{\mathbf i})=\begin{cases}
{\mathbf B}(\til\lm)&\text{ if }\lm[i_1,\cd,i_k]
\text{ is a Young diagram for any }k=1,2,\cd,N,\\
\emptyset&\text{otherwise.}
\end{cases}
\]

\begin{thm}[\cite{HK,N1}]
\label{LRrule}
Let $\lambda$ and $\mu$ be Young diagrams with at most $n$ rows.
Then we have
\begin{equation}
\mathbf{B}(\lambda) \otimes \mathbf{B}(\mu) \cong  
\bigoplus_{\tiny\begin{array}{l}T\in \mathbf{B}(\mu),\\
{\rm ME}(T)= \fsq(i_1) 
\otimes\cdots\otimes\fsq(i_N)
\end{array}} 
\mathbf{B}(\lambda:i_1,i_2,\ldots,i_N).
\label{LR}
\end{equation}
\end{thm}
Define
\[
\mathbf{B}(\mu)_\lm^\nu
:= 
\left\{\begin{array}{c|l}T\in{\mathbf B}(\mu)&
\begin{array}{l} 
{\rm ME}(T)=\MECry. \\
\text{For any }k=1,\cd,N,\\
\lm[i_1,\cd,i_k]\text{ is a Young diagram and}\\
\lm[i_1,\cd,i_N]=\nu.
\end{array}
\end{array}
\right\},
\]
which is called a {\it Littlewood-Richardson crystal}
associated with a triplet $(\lm,\mu,\nu)$.

\renewcommand{\thesection}{\arabic{section}}
\section{Robinson-Schensted-Knuth(RSK) correspondence }
\setcounter{equation}{0}
\renewcommand{\theequation}{\thesection.\arabic{equation}}

In this section we review the 
Robinson-Schensted-Knuth(RSK) correspondence with respect to 
column bumping procedure. For the contents of this 
section see \cite{F} (in particular, Appendix A.).
\subsection{Column Bumping and RSK Correspondence}

For an integer $x$ and a Young tableau $T$, we define the
column bumping procedure:
\begin{df}
\begin{enumerate}
\item
\begin{enumerate}
\item
If all entries in the 1-st column of $T$ are greater than $x$,
put $x$ just beneath the 1-st column and the procedure is over.
\item
Otherwise, let $y$ be the top entry in the 1-st column that
is equal to or smaller than $x$ and put $x$ in the box 
and bump the entry $y$ out.
\item
Do the same one for $y$ and the second column. If it does not
stop at the last column, make a new box next to the last column
and put the entry in the new box. 
\end{enumerate}
We denote the resulting tableau by $x\to T$. 
\item
The shape of $x\to T$ is a diagram added one box to the
original shape of $T$. We shall denote the added new 
box by New$(x)$ and call the new box by $x$.
\end{enumerate}
\end{df}
The following lemma is known as the 'column bumping lemma'
\begin{lem}\label{col-lem}
Let $T$ be a tableau and $x,x'$ positive integers. 
In the column bumping 
$x'\to(x\to T)$, we have:
\begin{enumerate}
\item
If $x<x'$, then New$(x')$ is weakly left of and 
strictly below New$(x)$.
\item
If $x\geq x'$, then New$(x)$ is strictly left of and 
weakly below New$(x')$.
\end{enumerate}
\end{lem}
It is shown similarly to the row bumping lemma(\cite{F}).

As is well-known that there is the reverse operation of this 
procedure, which is called an reverse (column) bumping.

\begin{df}
A two-rowed array
$w=\begin{pmatrix}u_1u_2\cd u_m\\ 
v_1v_2\cd v_m\end{pmatrix}$ is in a lexicographic order 
(of column type) if it satisfies:
(i) $u_1\leq u_2\leq \cd\leq u_m$.
(ii)
If $u_k=u_{k+1}$, then $v_k\geq v_{k+1}$.
\end{df}
Let $w$ be a two-rowed array in a lexicographic order with 
length $m$ as above.
We call the following procedure the RSK procedure:
\begin{enumerate}
\item
Set $P_1=v_1$ and $Q_1=u_1$.
\item
We obtain $(P_{k+1},Q_{k+1})$ from $(P_k,Q_k)$ by 
$P_{k+1}=v_{k+1}\to P_k$ and put $u_{k+1}$ to the same place 
in $Q_k$ as the newbox by $v_{k+1}$ in $P_{k+1}$.
\item
Set $R(w):=(P,Q)=(P_m,Q_m)$.
\end{enumerate}
Note that $P$ and $Q$ are Young tableaux with entries $1,\cd, m$
and same shape. We call a tableau $Q$ a recording tableau of $P$.
This procedure is reversible by using the
reverse column bumping:
For a pair of Young tableaux $(P,Q)$, we apply the reverse bumping to
$P$ starting from the box in $P$ which is in the same position as the box
with the right-most maximum entry in $Q$ and remove the entry from $Q$.
Repeat this procedure until the tableaux become empty.
We obtain the two-rowed array from $(P,Q)$, which gives the reverse
of the RSK procedure.
\begin{thm}(RSK correspondence)
Let ${\mathbf W}[n;m]$ be the set of two-rowed array in 
lexicographic order (of column type) with length $m$ and 
entries $1,\cd,n$ and ${\mathbf P}[n;m]$ be the set of 
pairs of same-shaped 
Young tableaux with $m$ boxes and entries $1,\cd,n$.
Then the map $R$ as above gives a bijection between
${\mathbf W}[n;m]$ and ${\mathbf P}[n;m]$.
\end{thm}

\subsection{Knuth equivalence and Crystal equivalence}

In this article, a {\it word} means a finite sequence of 
non-negative integers. 
\begin{df}(Knuth equivalence)
\begin{enumerate}
\item
Each of the following transformations between 3-letter words is
called a fundamental Knuth transformation:
\begin{enumerate}
\item[$\bullet$]
$K:yxz\longleftrightarrow yzx$ if $x<y\leq z$
\item[$\bullet$]
$K':xzy\longleftrightarrow zxy$ if $x\leq y<z$.
\end{enumerate}
\item
If two words with same length $w$ and $w'$ are Knuth 
equivalent if one can be transformed to the other 
by a sequence of the elementary Knuth transformations and 
we denote it by $w\douchik w'$.
\end{enumerate}

\end{df}
Here let us mention the relation between the 
crystal ${\mathbf B}$ and the Knuth equivalence.
The following lemma is well-known:
\begin{lem}
There exists the following non-trivial isomorphism 
of crystals:
${\mathbf R}:{\mathbf B}\ot {\mathbf B}\ot {\mathbf B}
\to {\mathbf B}\ot {\mathbf B}\ot {\mathbf B}$ by :
\begin{eqnarray*}
&&{\mathbf R}(\fsq(b) \otimes \fsq(a) \otimes \fsq(c) ) 
= \fsq(b) \otimes \fsq(c) \otimes \fsq(a),  \quad
 {\mathbf R}(\fsq(b) \otimes \fsq(c) \otimes \fsq(a)) 
= \fsq(b) \otimes \fsq(a) \otimes \fsq(c) \q
\text{ if }a\leq b<c,\\
&&{\mathbf R}(\fsq(c) \otimes \fsq(a) \otimes \fsq(b))
= \fsq(a) \otimes \fsq(c) \otimes \fsq(b), \quad
 {\mathbf R}(\fsq(a) \otimes \fsq(c) \otimes \fsq(b))
= \fsq(c) \otimes \fsq(a) \otimes \fsq(b)
\q \text{ if }a<b\leq c,\\
&&{\mathbf R}={\rm id}, \q\text{ otherwise.}
\end{eqnarray*}
\end{lem}
This is known as a combinatorial R matrix.
Indeed, 
\[
 {\mathbf B}^{\ot 3}\cong 
B(\hthreebox(,,))\oplus B\left(\twoone(,,)\,\,\right)^{\oplus 2}
\oplus B\left(\vthreebox(,,)\,\right),
\]
and the map ${\mathbf R}$ flips two components 
$B(\twoone(,,))$ each other.
Using this, we induce certain equivalent relation
between elements in ${\mathbf B}^{\ot m}$.
\begin{df}(Crystal equivalence)
Two elements $b,b'$ in ${\mathbf B}^{\ot m}$ are 
{\it crystal equivalent} if one is obtained by the others by 
applying a sequence of ${\mathbf R}$'s and we denote it by 
$b\, \douchic \, b'$.
\end{df}
The following is trivial by the theory of crystal bases:
\begin{pro}\label{c-e-o}
If $b\,\,\douchic\,\,b'$ ($b,b'\in {\mathbf B}^{\ot m}$), then
$\eit b\,\,\douchic\,\, \eit b'$ or $\eit b=\eit b'=0$
(resp. $fit b\,\,\douchic\,\,
 \fit b'$ or $\fit b=\fit b'=0$) for
any $i$.
\end{pro}
By the definitions we can easily see:
\begin{lem}\label{k-c}
For words $w=a_1a_2\cd a_m$ and $w'=b_1b_2\cd b_m$, set 
$b:=\fsq(a_m)\ot\cd\ot \fsq(a_1)$ and 
$b':=\fsq(b_m)\ot\cd\ot \fsq(b_1)$. Then we have 
$w\,\douchik\, w'$ if and only if $b\,\douchic \, b'$.
\end{lem}
\begin{df}
For a skew Young tableau $S$, a word $w(S)$ is defined by reading 
the entries in each row from left to right and 
from the bottom row to the top row, which is 
called a skew tableau word of $S$.
\end{df}
The following is given in \cite{F}.
\begin{pro}\label{douchi-w}
For a Young tableau $T$ and a positive integer $x$, we have
$ w(x\to T)\douchik x\cdot w(T), $
and then for positive integers $x_1,\cd,x_m$ we have
\[
 w(x_1\to(x_2\to(\cd(x_{m-1}\to x_m))))\,\,\douchik
\,\,x_1x_2\cd x_{m-1}x_m.
\]
\end{pro}

\renewcommand{\thesection}{\arabic{section}}
\section{Main Theorem}
\setcounter{equation}{0}
\renewcommand{\theequation}{\thesection.\arabic{equation}}

Let $\kappa^i$ ($i=1,2$) be skew diagrams such that
$|\kappa^1|=|\kappa^2|=:N$
and $\lm^i,\nu^i$ ($i=1,2$) be Young diagrams satisfying
$\kappa^i=\nu^i\setminus\lm^i$. Now, 
let us define the the map ${\mathcal S}$:
\[
{\mathcal S} : \mathbf{P}(\kappa^1,\kappa^2) \to
 \coprod_{\mu} \left(\mathbf{B}(\mu)^{\nu^1}_{\lambda^1}
 \times \mathbf{B}(\mu)^{\nu^2}_{\lambda^2}\right)
\qquad (f \; \mapsto \; (T^1,T^2)),
\]
where $\mu$ runs over the set of Young diagrams with 
$|\kappa^1|=|\kappa^2|=|\mu|(=N)$.

Set
\begin{eqnarray*}
&&{\mathbf S}(\kappa^1,\kappa^2):=\left\{
\begin{array}{l|l}
S&
\begin{array}{l}
S\text{ is a skew tableau of shape }\kappa^1\text
{ and the number of entry }i\text{ is }\kappa^2_i,\\
ME(S)=\MECry\text{ satisfies that }\lm^2[i_1,\cd,i_k]
\text{ is}\\
\text{a Young diagram for $k=1,\cd,N$ 
and }\lm^2[i_1,\cd,i_N]=\nu^2.
\end{array}
\end{array}
\right\},\\
&&{\mathbf W}(\kappa^1,\kappa^2):=\left\{
\begin{array}{l|l}
  w = \left(
 \begin{array}{c}
w^1\\w^2
\end{array} \right)&
\begin{array}{l}
w\text{ is a lexicographic two-rowed array of length }N,\\
\sharp\{i\in w^j\}=\kappa^j_i\,\,(j=1,2),\\
\text{the column bumping of }w^2\text{ is in }
{\bf B}(\mu)_{\lm^2}^{\nu^2} \text{ and }\\
\text{the recording tableau by $w^1$ is in }
{\bf B}(\mu)_{\lm^1}^{\nu^1}.
\end{array}
\end{array}
\right\}
\end{eqnarray*}
where an element in ${\mathbf S}(\kappa^1,\kappa^2)$ is called 
a Littlewood-Richardson skew tableau associated with 
$(\kappa^1,\kappa^2)$.
Let us define maps:
\[
\hspace{-20pt}
{\mathcal S}_1:{\bf P}(\kappa^1,\kappa^2)\to 
{\bf S}(\kappa^1,\kappa^2),\q
{\mathcal S}_2:{\bf S}(\kappa^1,\kappa^2)\to 
{\bf W}(\kappa^1,\kappa^2),\q
{\mathcal S}_3:{\bf W}(\kappa^1,\kappa^2)\to 
\coprod_{\mu} \left( \mathbf{B}(\mu)^{\nu^1}_{\lambda^1}
 \times \mathbf{B}(\mu)^{\nu^2}_{\lambda^2}\right).
\]
\begin{df}
\label{def-S}
\begin{enumerate}
\item
For a picture $f=(f_1,f_2)\in \mathbf{P}(\kappa^1,\kappa^2)$
(where $f_1,f_2$ means a coordinate of a box in $\kappa^2$),
let $S$ be a skew tableau of shape $\kappa^1$ 
whose $(i,j)$-entry 
$S_{i,j}=f_1(i,j)$. Define ${\mathcal S}_1(f):=S$.
\item
For $S\in{\mathbf S}(\kappa^1,\kappa^2)$, 
writing $ME(S)=\fsq(a_1) \otimes \fsq(a_2) \otimes 
\cdots \otimes \fsquare(4mm,a_N)$, define a word 
$w^2=a_1a_2\cd a_N$. 
Let $b_i$ ($i=1,2,\cd,N$) be the row number of the place of 
$a_i$ in $S$ and set $w^1=b_1b_2\cd b_N$.
Define 
\[
 {\mathcal S}_2(S):=w = \left(\begin{array}{c}
w^1\\w^2
\end{array} \right) =
\left(
\begin{array}{cccc}
b_1 & b_2 & \cdots & b_N\\
a_1 & a_2 & \cdots & a_N\\
\end{array}
\right).
\]
\item
For the two rowed array 
$w = \left(\begin{array}{c}
w^1\\w^2
\end{array} \right) =
\left(
\begin{array}{cccc}
b_1 & b_2 & \cdots & b_N\\
a_1 & a_2 & \cdots & a_N\\
\end{array}
\right)\in {\mathbf W}(\kappa^1,\kappa^2)$,
apply the column bumping procedure to $w^2$ and obtain 
a tableau $T^2=a_N\to(\cd(a_2\to a_1))$. 
Let $T^1$ be the recording tableau 
of $T^2$ using $w^1$. Define ${\mathcal S}_3(w)=(T^1,T^2)$.
\item
Finally, define ${\mathcal S}={\mathcal S}_3\circ
{\mathcal S}_2\circ {\mathcal S}_1$.
\end{enumerate}
\end{df}

Next, let us define a map ${\mathcal C}$
\[
 {\mathcal C}:\coprod_{\mu} \left(
 \mathbf{B}(\mu)^{\nu^1}_{\lambda^1}
 \times \mathbf{B}(\mu)^{\nu^2}_{\lambda^2}\right)\to
{\bf P}(\kappa^1,\kappa^2).
\]
To carry out this task, we define the following maps:
\[
 \hspace{-20pt}{\mathcal C}_3:\coprod_{\mu} \left(
 \mathbf{B}(\mu)^{\nu^1}_{\lambda^1}
 \times \mathbf{B}(\mu)^{\nu^2}_{\lambda^2}\right)
\to 
{\bf W}(\kappa^1,\kappa^2),\q
{\mathcal C}_2:{\bf W}(\kappa^1,\kappa^2)\to 
{\bf S}(\kappa^1,\kappa^2),\q
{\mathcal C}_1:{\mathbf S}(\kappa^1,\kappa^2)\to 
{\bf P}(\kappa^1,\kappa^2).
\]
\begin{df}
\begin{enumerate}
\item
For a pair of tableaux $(T^1,T^2)\in 
\coprod_{\mu\left(\mathbf{B}(\mu)^{\nu^1}_{\lambda^1}\times
\mathbf{B}(\mu)^{\nu^2}_{\lambda^2}\right)}$, apply the reverse column
bumping to $T^2$ by using $T^1$ as a recording tableau and 
set $c_N c_{N-1}\cd c_1$ a sequence obtained from $T^2$
($c_i$ is the $N+1-i$-th entry bumped out from $T^2$.).
Set $w^2:=c_1,\cd,c_N$ and let $d_i$ be 
the entry in the same place in $T^1$ as the $(N-i+1)$-th 
removed box in $T^2$
and set $w^1:=d_1,\cd,d_N$. 
Define ${\mathcal C}_3(T^1,T^2)=w
=\left(\begin{array}{c}w^1\\w^2\end{array}\right)$.
\item
For 
\[
 w=\left(\begin{array}{c}w^1\\w^2\end{array}\right)=
\begin{pmatrix}d_1d_2\cd d_N\\ 
c_1c_2\cd c_N\end{pmatrix}\in {\mathbf W}(\kappa^1,\kappa^2),
\]
put $c_1c_2\cd c_N$ to $\kappa^1$ according 
to the middle-eastern
ordering and set $S$ the resulting skew tableau, whose shape is 
$\kappa^1$. Define ${\mathcal C}_2(w)=S$.
\item
For $S\in{\mathbf S}(\kappa^1,\kappa^2)$, 
define ${\mathcal C}_1(S)=f$ by 
$f(i,j):=(S_{ij},\lambda^2_{S_{ij}}+p(S;i,j))$ for 
$(i,j)\in \kappa^1$, where 
$p(S;i,j)$ is as above and $S_{ij}$ is the $(i,j)$-entry of $S$.
\item
Finally, we define $\mathcal C=\mathcal C_1\circ\mathcal
 C_2\circ\mathcal C_3$.
\end{enumerate}
\end{df}
Note that well-definedness of each map will be shown later.
\begin{thm}
\label{main}
In the above setting, the maps ${\mathcal S}$ and $\mathcal C$ 
are both well-defined bijective maps between 
$\mathbf{P}(\kappa^1,\kappa^2)$ and 
$\coprod_{\mu} \left( \mathbf{B}(\mu)^{\nu^1}_{\lambda^1}
\times \mathbf{B}(\mu)^{\nu^2}_{\lambda^2} \right)$
and  inverse
each other.
\end{thm}
Here note that the set $\coprod_{\mu} 
\left( \mathbf{B}(\mu)^{\nu^1}_{\lambda^1}
\times \mathbf{B}(\mu)^{\nu^2}_{\lambda^2} \right)$
consists of pairs of same shaped Young tableaux, which means that 
this theorem is an analogue of the RSK correspondence.
\begin{ex}
We take the following skew diagrams:

\setlength{\unitlength}{1pt}
\begin{picture}(100,40)

\put(0,15){{$\kappa^1 =$}}
\put(30,0){\framebox(30,10){}}
\put(40,0){\line(0,1){20}}
\put(50,0){\line(0,1){30}}
\put(60,10){\line(0,1){20}}
\put(70,20){\line(0,1){10}}
\put(40,20){\line(1,0){30}}
\put(50,30){\line(1,0){20}}

\put(90,15){{$\kappa^2 =$}}
\put(120,-5){\framebox(20,20){}}
\put(130,-5){\line(0,1){20}}
\put(120,5){\line(1,0){20}}
\put(140,15){\framebox(10,20){}}
\put(140,25){\line(1,0){20}}
\put(150,25){\line(1,0){10}}
\put(160,25){\line(0,1){10}}
\put(150,35){\line(1,0){10}}
\end{picture}

Let $f_a\in{\mathbf P}(\kappa^1,\kappa^2)$ be \\
\qq\qq$\qq f_a =$
\begin{tabular}{c||c|c|c|c|c|c|c}
$\kappa^1$ & $(1,3)$ & $(1,4)$ & $(2,2)$ & $(2,3)$ & $(3,1)$ & $(3,2)$ & $(3,3)$\\
\hline
$\kappa^2$ & $(1,3)$ & $(3,1)$ & $(1,4)$ & $(3,2)$ & $(2,3)$ & $(4,2)$ & $(4,1)$\\
\end{tabular}\\
Here we have 

\begin{picture}(80,40)

\put(-10,10){$S_a={\mathcal S}_1(f_a) =$}
\put(55,1){$\hthreebox(2,4,4)$}
\put(63.5,9.5){$\htwobox(1,3)$}
\put(72,18){$\htwobox(1,3)$}
\put(100,10){\text{and then ${\rm ME}(S_a)=\fsq(3) \otimes 
\fsq(1) \otimes \fsq(3) \otimes \fsq(1) \otimes 
\fsq(4) \otimes \fsq(4) \otimes \fsq(2)$.}}
\end{picture}

\nd
Then we get 
$w_a={\mathcal S}_2(S_a)=\begin{pmatrix}1122333\\ 3131442
\end{pmatrix}$ and then finally, we have

\[T^2 : \fsq(3) \dashrightarrow \htwobox(1,3) 
\dashrightarrow \twoone(1,3,3)
\dashrightarrow  \threeone(1,1,3,3) 
\dashrightarrow \threeoneone(1,1,3,3,4)
\dashrightarrow \threetwoone(1,1,3,3,4,4) 
\dashrightarrow \threethreeone(1,1,3,2,3,4,4) = T_a^2\qq\qq\]
\[T^1 : \fsq(1) \dashrightarrow \htwobox(1,1) 
\dashrightarrow \twoone(1,1,2)
\dashrightarrow  \threeone(1,1,2,2) 
\dashrightarrow \threeoneone(1,1,2,2,3)
\dashrightarrow \threetwoone(1,1,2,2,3,3) 
\dashrightarrow \threethreeone(1,1,2,2,3,3,3) = T_a^1,\qq\qq\]
that is, ${\mathcal S}_3(w_a)=(T^1_a,T^2_a)$.\\
Conversely, For $(T^1,T^2)=
\left(\threethreeone(1,1,2,2,3,3,3), 
\threethreeone(1,1,3,2,3,4,4)\right)$, applying 
the reverse column bumping to $T^2$ using $T^1$, we get
$c_7=2,\,c_6=4,\,c_5=4,\,c_4=1,\,c_3=3,\,c_2=1,\,c_1=3$ and 
$d_1=d_2=1,\,d_3=d_4=2,d_5=d_6=d_7=3$ and then 
\[
w= {\mathcal C}_3(T^1,T^2)=\begin{pmatrix}1122333\\ 3131442
\end{pmatrix}.
\]
We obtain 

\nd
\begin{picture}(80,40)

\put(0,15){$S ={\mathcal C}_2(w)=$}
\put(60,5){$\hthreebox(c_7,c_6,c_5)$}
\put(68.5,13.5){$\htwobox(c_4,c_3)$}
\put(77,22){$\htwobox(c_2,c_1)$}

\put(110,10){$=$}

\put(130,5){$\hthreebox(2,4,4)$}
\put(138.5,13.5){$\htwobox(1,3)$}
\put(147,22){$\htwobox(1,3)$}
\put(180, 13.5){\text{ and then finally, we have  }}
\end{picture}
\qq\qq\qq \qq\qq\qq

\nd 
${\mathcal C}_1(S) = $
\begin{tabular}{c||c|c|c|c|c|c|c}
$\kappa^1$ & $(1,3)$ & $(1,4)$ & $(2,2)$ & $(2,3)$ & $(3,1)$ & $(3,2)$ & $(3,3)$\\
\hline
$\kappa^2$ & $(1,3)$ & $(3,1)$ & $(1,4)$ & $(3,2)$ & $(2,3)$ & $(4,2)$ & $(4,1)$\\
\end{tabular} $ = f_a$.

\end{ex}
\medskip

To show the theorem, it suffices to prove 
\begin{enumerate}
\item
The well-definedness of ${\mathcal S}$.
\item
The well-definedness of ${\mathcal C}$.
\item
Bijectivity of ${\mathcal S}$ and ${\mathcal C}$.
\end{enumerate}
We shall show these in the subsequent sections.

\renewcommand{\thesection}{\arabic{section}}
\section{Well-definedness of $\mathcal S$}
\setcounter{equation}{0}
\renewcommand{\theequation}{\thesection.\arabic{equation}}

For the well-definedness of $\mathcal S$,
we shall prove the following:
\begin{pro}\label{prop:S}
The maps ${\mathcal S}_i$ $(i=1,2,3)$ are well-defined.
\end{pro}
Indeed, the well-definedness of $\mathcal S_3$ is obvious by the
definition.
\subsection{Well-definedness of ${\mathcal S}_1$}

For $f\in {\mathbf S}(\kappa^1,\kappa^2)$, 
by the similar argument in \cite{NS,NS2}, 
we can show that $S:={\mathcal S}_1(f)$
is a skew tableau. Thus, we may show:
\begin{lem}
\label{ek-S}
For any $k=1\cd,n$ and the skew tableau $S={\mathcal S}_1(f)$, we have
\[
 \tilde{e}_k(ME(Y_{\lambda^2}) \otimes ME(S))=0,
\]
where $Y_{\lm^2}$ is a Young tableau of shape $\lm^2$ satisfying that 
all the entries in $k$-th row are $k$ $(k=1,\cd,n)$, which is called 
a highest tableau.
\end{lem}
{\sl Proof.}
Write
\[
 ME(Y_{\lambda^2}) \otimes ME(S)=\fsq(i_1)\ot\cd\ot\fsq(i_N).
\]
By the rule of the action of $\til e_k$, we may show 
\begin{equation}
 \sharp\{j|i_j=k,\,\,j\leq p\}
\geq  \sharp\{j|i_j=k+1,\,\,j\leq p\}
\label{number}
\end{equation}
for any $p=1,\cd,N$. In the skew diagram $\kappa^2$, we have

\begin{picture}(200,50)
\put(30,0){\framebox(40,15){$C$}}
\put(70,0){\framebox(50,15){$B$}}
\put(70,15){\framebox(50,15){$A$}}
\put(120,15){\framebox(30,15){$D$}}
\put(30,15){$\overbrace{\quad \quad \quad \quad}^{\lambda^2_k-\lambda^2_{k+1}}$}
\put(160,20){$\leftarrow k$-th row}
\put(160,5){$\leftarrow k+1$-th row}

\put(240,10){\Large{(in $\kappa^2$)}}
\end{picture}

For boxes $(k,j), (k+1,j)\in\kappa^2$, by the fact
$(k,j)\leqslant_P (k+1,j)$, we have
\[
 (x_1,y_1):=f^{-1}(k,j)\leqslant_J f^{-1}(k+1,j)=:(x_2,y_2).
\]
It is evident from the definition of the map $\mathcal S_1$ that 
\[
 S_{x_1,y_1}=k, \qq S_{x_2,y_2}=k+1.
\]
This implies that 
in the tensor product 
$ ME(Y_{\lambda^2}) \otimes ME(S)=\fsq(i_1)\ot\cd\ot\fsq(i_N)$, 
$k$'s from $A$ appear earlier than $k+1$'s from $B$ and then 
they are cancelled each other with respect 
to the action of $\til e_k$. 
In $ME(Y_{\lm^2})$, the number of $k$ exceeds the one of $k+1$
by $\lm^2_k-\lm^2_{k+1}$. Thus, $k+1$'s from the part 
$C$ in the figure also have been  
cancelled by $k$'s in $ME(Y_{\lm^2})$. 
Hence we obtain (\ref{number}) and then 
$ \tilde{e}_k(ME(Y_{\lambda^2}) \otimes ME(S))=0$
for any $k$.\qed

\nd
Thus, we have the well-definedness of $\mathcal S_1$.

\subsection{Well-definedness of ${\mathcal S}_2$}

First, let us show that the two-rowed array $w:={\mathcal S}_2(S)$ 
$(S\in{\mathbf S}(\kappa^1,\kappa^2))$ is in the
lexicographic order, that is, $b_1\leq b_\leq \cd\leq b_N$ and 
$a_j\geq a_{j+1}$ if $b_j=b_{j+1}$, where $a_j,b_j$ are as in 
Definition \ref{def-S}.
It follows immediately from the definition of $b_i$'s that 
$b_1\leq b_\leq \cd\leq b_N$. 
Let $k$ satisfy $b_1\leq k\leq b_N$ and 
$\{b_i,b_{i+1},\cd,b_{i+r}\}$ the maximal subsequence of $w^1$
such that $b_i=\cd=b_{i+r}=k$. This implies that 
$a_i,a_{i+1},\cd,a_{i+r}$ are the entries in $k$-th row of $S$. 
Since $S$ is a skew tableau, we obtain that 
$a_i\geq a_{i+1}\geq\cd\geq a_{i+r}$, which means that 
$w$ is in the lexicographic order.
Let $T^2$ be the tableau from $w^2$ by the column bumping. 
Let us show that $T^2\in \mathbf{B}(\mu)^{\nu^2}_{\lambda^2}$, {\it i.e.,}
\[
\tilde{e}_k(ME(Y_{\lambda^2}) \otimes ME(T^2))=0
\]
for any $k=1\cd,n$.
For this purpose, we see the following lemma.
\begin{lem}
\label{c-eq}
$ME(S)$ is crystal equivalent to $ME(T^2)$.
\end{lem}
{\sl Proof.}
For $w^2=a_1a_2\cd a_N$, since 
$T^2$ is obtained by the column bumping
procedure of $a_N\cd a_1$, we know that 
$w(S)=a_N a_{N-1} \cdots a_1\; \douchik \; w(T^2)$, which means
$ME(S)\douchic ME(T^2)$ by Lemma \ref{k-c}. \qed

By the Lemma \ref{c-eq}, we have $ME(S)\douchic ME(T^2)$ and 
then $ME(Y_{\lm^2})\ot ME(S)\douchic ME(Y_{\lm^2})\ot ME(T^2)$.
We also have 
\[
 \til e_k(ME(Y_{\lm^2})\ot ME(S))=0,
\]
for any $k$ by Lemma \ref{ek-S}. This and Proposition 
\ref{c-e-o} show that 
\[
 \til e_k(ME(Y_{\lm^2})\ot ME(T^2))=0,
\]
for any $k$ and then we have $T^2\in{\bf B}(\mu)^{\nu^2}_{\lm^2}$.

For $w:={\mathcal S}_2(S)$, 
we set $(T^1,T^2):={\mathcal S}_3(w)$. For our purpose, it suffices 
to show $T^1\in {\bf B}(\mu)^{\nu^1}_{\lm^1}$, 
that is, 
$ \til e_k(ME(Y_{\lm^2})\ot ME(T^1))=0$
for any $k$.

\begin{lem}
\label{cclem}
For $1\leq c_1,\cd,c_k\leq n$, assume that 
\[
b_1:=\cdots \otimes \Longsquare(0.4cm,0.8cm,c_{i-1}) \otimes
\fsquare(0.4cm,c_i) \otimes \Longsquare(0.4cm,0.8cm,c_{i+1}) \otimes 
\Longsquare(0.4cm,0.8cm,c_{i+2}) \otimes 
\cdots 
\; \douchic \;
\cdots \otimes \Longsquare(0.4cm,0.8cm,c_{i-1}) \otimes
 \Longsquare(0.4cm,0.8cm,c_{i+1})
\otimes \fsquare(0.4cm,c_i)  \otimes 
\Longsquare(0.4cm,0.8cm,c_{i+2}) \otimes 
\cdots =:b_2.
\]
Applying the column bumping procedure to both
$b_1$ and $b_2$, the place of the new box New$(c_i)$ (resp. 
New$(c_{i+1})$) from $b_1$
coincides with the one of the new box New$(c_{i+1})$
(resp. New$(c_i)$) from $b_2$.
\end{lem}
{\sl Proof.}
Set $x:=c_{i}$, $y:=c_{i-1}$ and $z:=c_{i+1}$.
First we consider the case $x\leq y<z$.
Let $T_p$ (resp. $T_q$)be the tableau obtained from $b_1$ 
(resp. $b_2$) by the column bumping procedure.
It follows immediately from the condition $x\leq y<z$ that 
\[
 w(T_p)\,\,\douchik \,\,c_k\cd zxy \cd c_1\,\,\douchik\,\,
c_k\cd xzy \cd c_1\,\,\douchik \,\,w(T_q),
\]
which shows that $T_p=T_q$.
Define the tableau $T'$ by column bumping 
\begin{eqnarray}
 T'&:=&z\to (x\to (y\to(\cd (c_2\to c_1))))
\label{11}\\
&=&x\to (z\to (y\to(\cd (c_2\to c_1)))).
\label{22}
\end{eqnarray}

Let $X={\rm New}(x)$ and $Z={\rm New}(z)$ be the new boxes
in each column bumping.
By the condition $x<z$, in the bumping (\ref{11}) we have
by the column bumping lemma:

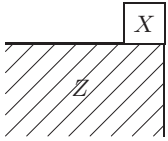
\begin{picture}(100,70)
\put(50,40){\framebox(15,15){$X$}}
\put(65,40){\line(0,-1){35}}
\put(50,40){\line(-1,0){45}}
\put(30,20){$Z$}
\hatchpattern{8}{\path}
    \hatch(5,5)(65,40)
\end{picture}

In (\ref{22}), we have

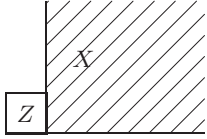
\begin{picture}(100,70)
\put(20,10){\framebox(15,15){$Z$}}
\put(35,10){\line(1,0){60}}
\put(35,25){\line(0,1){35}}
\put(45,35){$X$}
\hatchpattern{8}{\path}
    \hatch(35,10)(95,60)
\end{picture}

These mean that $X$ (resp. $Z$) in (\ref{11}) coincides
with $X$ (resp. $Z$) in (\ref{22}).
We can show the case $x<y\leq z$ and the case 
$x=c_i$, $z=c_{i+1}$ and $y=c_{i+2}$ 
similarly.\qed

To show $\tilde{e}_k(ME(Y_{\lambda^1}) \otimes ME(T^1))=0$
for any $k$, we see the $k$-th and $k+1$-th rows of $S$.

\vspace{-0.5cm}
\begin{picture}(200,50)
\put(30,0){\framebox(40,15){$d_1 \cdots d_j$}}
\put(70,0){\framebox(80,15){$\cdots$}}
\put(70,15){\framebox(80,15){$\cdots$}}
\put(150,15){\framebox(50,15){$a_{m+1} \cdots a_n$}}
\put(73,5){$b_1$}
\put(73,20){$a_1$}
\put(138,5){$b_m$}
\put(138,20){$a_m$}
\put(210,20){$\leftarrow k$-th row}
\put(210,5){$\leftarrow k+1$-th row}
\put(300,10){\Large{(in $S$)}}
\end{picture}

By this figure, we know that 
\[
 a_1< b_{i-1}\leq b_i\q i=2,3,\cd, m.
\]
This induces the following transformations of $ME(S)$
by the map $\mathbf R$:
\begin{eqnarray*}
ME(S)&=& \cdots \otimes \fsquare(0.5cm,a_2) \otimes
\fsquare(0.5cm,a_1) \otimes \fsquare(0.5cm,b_m) 
\otimes \Longsquare(0.5cm,0.8cm,b_{m-1}) \otimes
\cdots \otimes \fsquare(0.5cm,b_1) \otimes \fsquare(0.5cm,d_j) 
 \otimes \cdots\\
&\douchic&
\cdots \otimes \fsquare(0.5cm,a_2) \otimes
\fsquare(0.5cm,b_m) \otimes \fsquare(0.5cm,a_1) 
\otimes \Longsquare(0.5cm,0.8cm,b_{m-1}) \otimes
\cdots \otimes \fsquare(0.5cm,b_1) \otimes \fsquare(0.5cm,d_j) 
 \otimes \cdots\\ &&\qq \cd\cd\cd\cd \\
&\douchic&
\cdots \otimes \fsquare(0.5cm,a_2) \otimes
\fsquare(0.5cm,b_m) \otimes \Longsquare(0.5cm,0.8cm,b_{m-1}) 
\otimes\cd\otimes \fsquare(0.5cm,b_2) \ot
\fsquare(0.5cm,a_1)\ot \fsquare(0.5cm,b_1)\ot
 \fsquare(0.5cm,d_j) 
 \otimes \cdots.
\end{eqnarray*}
Furthermore, we have $a_j<b_{i-1}\leq b_i$ for 
$2\leq j< i\leq m$. Thus, 
repeating the above transformations we get 
\begin{equation}
ME(S)\douchic
\cdots \otimes\fsquare(0.5cm,{a_m}) 
\otimes \fsquare(0.5cm,{b_m})
\otimes \Longsquare(0.5cm,0.8cm,{a_{m-1}}) \otimes
\Longsquare(0.5cm,0.8cm,{b_{m-1}})
\otimes \cdots \otimes
\fsquare(0.5cm,{a_2}) \otimes
\fsquare(0.5cm,{b_2}) \otimes
\fsquare(0.5cm,{a_1}) \otimes
\fsquare(0.5cm,{b_1})
\otimes \fsquare(0.5cm,d_j) \otimes \cdots
=:w',
\label{last-c}
\end{equation}
which means that the resulting tableaux 
 by column bumping of $ME(S)$ and $w'$ 
are same as $T^2$ by Lemma \ref{cclem}.
Considering the column bumping of $w'$, 
set $A_1:={\rm New}(a_1)$
and $B_1:={\rm New}(b_1)$ in $T^2$.
We have 

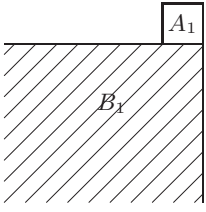
\setlength{\unitlength}{1pt}
\begin{picture}(90,90)
\put(60,60){\framebox(15,15){$A_1$}}
\put(60,60){\line(-1,0){60}}
\put(75,60){\line(0,-1){60}}
\put(35,35){$B_1$}
\hatchpattern{8}{\path}
    \hatch(0,0)(75,60)
\end{picture}

Since the entry $a_1$ (resp. $b_1$) has been placed at
the $k$ (resp. $k+1$)-th row in $S$, 
in $T^1$ we have

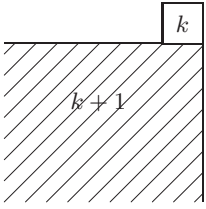
\setlength{\unitlength}{1pt}
\begin{picture}(90,90)
\put(60,60){\framebox(15,15){$k$}}
\put(60,60){\line(-1,0){60}}
\put(75,60){\line(0,-1){60}}
\put(25,35){$k+1$}
\hatchpattern{8}{\path}
    \hatch(0,0)(75,60)
\end{picture}

So, in $ME(T^1)$ the $k$ as above appears earlier than the $k+1$.
We know that the positions of 
 $New(a_i)$ and $New(b_i)$ in $T^1$ are in the similar relation
to the one of New$(a_1)$ and New$(b_1)$
and then in $ME(T^1)$ the $k$'s from $a_1\cd, a_m$  cancel
$k+1$'s from $b_1\cd,b_m$.
Moreover, in $ME(Y_{\lm^1})$ there are 
more $\lm^1_k-\lm^1_{k+1}$ $k$'s 
than $k+1$'s. Thus, $k+1$'s from $d_1\cd,d_j$ have been 
cancelled in $ME(T^1)$ and this implies 
$\til e_k(ME(Y_{\lm^2})\ot ME(T^1))=0$ for any $k$.
Now, we obtain $T^1\in \mathbf{B}(\mu)^{\nu^1}_{\lambda^1}$
and the well-definedness of the map $\mathcal S_2$ and then 
$\mathcal S$, which
 completes the proof of Proposition~\ref{prop:S}.
\qed

\renewcommand{\thesection}{\arabic{section}}
\section{Well-definedness of $\mathcal C$}
\setcounter{equation}{0}
\renewcommand{\theequation}{\thesection.\arabic{equation}}

To show the well-definedness of the map $\mathcal C$, we should
prove that $f:={\mathcal C}(T^1,T^2)$ is a PJ-picture
from $\kappa^1$ to $\kappa^2$.
In the course of the proof, we shall also show that 
the maps $\mathcal C_1$, $\mathcal C_2$ and $\mathcal C_3$
are well-defined. Indeed, that of $\mathcal C_3$ is immediate 
from the definition. 
\begin{pro}\label{prop:C}
Let $S$ be the filling of shape $\kappa^1$
appearing in the definition of 
$\mathcal C_2$. Then $S$ is a skew tableau of shape $\kappa^1$.
\end{pro}
{\sl Proof.}
For $w=\begin{pmatrix}w^1\\w^2
\end{pmatrix}\in{\mathbf W}(\kappa^1,\kappa^2)$, set
$(T^1,T^2):={\mathcal S}_3(w)$, which is in 
$\coprod_{\mu} 
\left( \mathbf{B}(\mu)^{\nu^1}_{\lambda^1}
\times \mathbf{B}(\mu)^{\nu^2}_{\lambda^2} \right)$
as we have seen in the last section.
Since $T^1$ is in ${\bf B}(\mu)_{\lm^1}^{\nu^1}$, the number of 
entry $k$'s$(k=1,\cd,n)$ is $h:=\nu^1_k-\lm^1_k$. 
Let $X_1\cd,X_h$ be the positions of $k$'s 
in $T^1$ from right to 
left. Note that ${(T^1)}^{(k)}=\{X_1,\cd,X_h\}$.
And let $x_j$ ($j=1\cd,h$) be the entry in $T^2$ at the same 
position as $X_j$.
By the definition of $\mathcal C_2$, 
the entries in $k$-th row of 
$S$ consist of elements obtained by 
reverse column bumping, that is,
the entry $S_{k,\lm^1+i}$ is the element
by the inverse column bumping of $x_i$.

Now, assume that $S_{k,\lm^1+i}>S_{k,\lm^1+i+1}$.
In the column bumping of $w^2=ME(S)$ to $T^2$,
the new box by $S_{k,\lm^1+i}$ (resp. $S_{k,\lm^1+i+1}$) 
has $x_i$ (resp. $x_{i+1}$) and it is placed
at $X_i$ (resp. $X_{i+1}$).
Applying the column bumping lemma (Lemma \ref{col-lem})
to these new boxes, we have

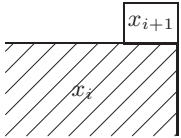
\begin{picture}(100,70)

\put(50,40){\framebox(20,15){$x_{i+1}$}}
\put(70,40){\line(0,-1){35}}
\put(50,40){\line(-1,0){45}}
\put(30,20){$x_i$}
\hatchpattern{8}{\path}
    \hatch(5,5)(70,40)
\end{picture}

This contradicts to the fact that 
$x_i$ is on the right side of $x_{i+1}$ and shows that 
$S_{k,\lm^1+i}\leq S_{k,\lm^1+i+1}$.

Next, let us check the condition for vertical directions in $S$.
Suppose that $S_{k,j}\geq S_{k+1,j}$. 
Then in $S$ we obtain the following $A,B$:

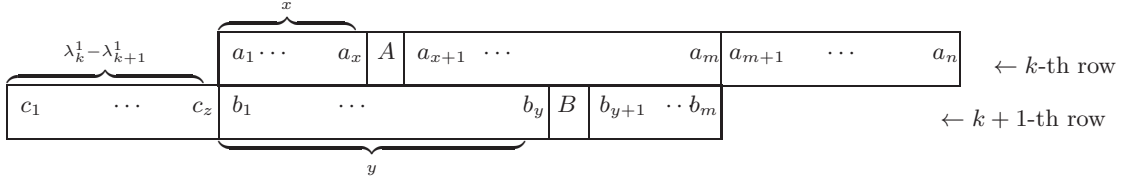
\setlength{\unitlength}{1pt}
\begin{picture}(400,80)(0,0)
\put(80,50){\framebox(280,20){}}
\put(0,30){\framebox(270,20){}}
\put(80,50){\line(0,-1){20}}
\put(270,30){\line(0,1){40}}
\put(85,60){$a_1$}
\put(95,60){$\cdots$}
\put(125,60){$a_x$}
\put(136,50){\line(0,1){20}}
\put(150,50){\line(0,1){20}}
\put(140,60){$A$}
\put(155,60){$a_{x+1}$}
\put(180,60){$\cdots$}
\put(258,60){$a_m$}
\put(273,60){$a_{m+1}$}
\put(310,60){$\cdots$}
\put(350,60){$a_n$}
\put(5,40){$c_1$}
\put(40,40){$\cdots$}
\put(70,40){$c_z$}
\put(85,40){$b_1$}
\put(125,40){$\cdots$}
\put(195,40){$b_y$}
\put(205,30){\line(0,1){20}}
\put(220,30){\line(0,1){20}}
\put(208,40){$B$}
\put(224,40){$b_{y+1}$}
\put(249,40){$\cdots$}
\put(258,40){$b_m$}

\put(80,70){$\overbrace{ \;\; \quad \quad \quad \quad  \quad}^{x}$}
\put(0,50){$\overbrace{\quad \quad \quad \qquad \qquad \quad}^{\lambda^1_k - \lambda^1_{k+1}}$}
\put(80,30){$\underbrace{\qquad \qquad \qquad \qquad \qquad \qquad \; \:}_{y}$}

\put(350,35){ $ \leftarrow k+1$-th row}
\put(370,55){ $ \leftarrow k$-th row}

\end{picture}

such that $A\geq B$, 
$a_i<b_j$ for $i\leq j$, $i=1,\cd,x$ and $j=1,\cd,m$.
Indeed, we get these by the following way.
\begin{enumerate}
\item 
Find the left-most pair ($a_s,b_s$) such that $a_s\geq b_s$.
\item
If $a_s\geq b_m$, then set $A:=a_s$ and $B:=b_m$.
\item
Otherwise,  compare $a_s$ and $b_{m-1}$ and if $a_s\geq b_{m-1}$, 
then set $A:=a_s$ and $B:=b_{m-1}$. 
\item
Otherwise, repeat the above procedure until getting 
$a_s\geq b_j$ for $j\geq s$. Then set $A:=a_s$ and $B:=b_j$.
\end{enumerate}
Since we have $a_1<b_{j-1}\leq b_j$ for $j=2,\cd,m$, and 
$a_1<B\leq b_{y+1}$ we have
\begin{eqnarray*}
ME(S) = \cdots \otimes \fsquare(0.5cm,a_n) \otimes \cdots 
\cdots &\otimes& \Longsquare(0.5cm,0.8cm,a_{x+1})  \otimes \fsquare(0.5cm,A) \otimes 
\fsquare(0.5cm,a_x) \otimes \cdots \otimes 
\fsquare(0.5cm,a_1)
\otimes \fsquare(0.5cm,b_m) \otimes 
\Longsquare(0.5cm,0.8cm,b_{m-1}) \otimes \cdots\\
 & \otimes&  \Longsquare(0.5cm,0.8cm,b_{y+1}) \otimes
\fsquare(0.5cm,B) \otimes \fsquare(0.5cm,b_y) 
\cdots \otimes \fsquare(0.5cm,b_1)
\otimes \fsquare(0.5cm,c_z) \otimes \cdots 
\otimes \fsquare(0.5cm,c_1) \otimes \cdots\\
\douchic \cdots \otimes \fsquare(0.5cm,a_n) \otimes \cdots 
\cdots& \otimes& \Longsquare(0.5cm,0.8cm,a_{x+1})  \otimes \fsquare(0.5cm,A) \otimes 
\fsquare(0.5cm,a_x) \otimes \cdots \otimes 
\fsquare(0.5cm,{b_m}) 
\otimes \fsquare(0.5cm,{a_1}) \otimes 
\Longsquare(0.5cm,0.8cm,b_{m-1}) \otimes \cdots\\
&  \otimes&  \Longsquare(0.5cm,0.8cm,b_{y+1}) \otimes
\fsquare(0.5cm,B) \otimes \fsquare(0.5cm,b_y) 
\cdots \otimes \fsquare(0.5cm,b_1)
\otimes \fsquare(0.5cm,c_z) \otimes \cdots 
\otimes \fsquare(0.5cm,c_1) \otimes \cdots\\
\douchic \cdots \otimes \fsquare(0.5cm,a_n) \otimes \cdots \otimes
\cdots& \otimes& \Longsquare(0.5cm,0.8cm,a_{x+1}) \otimes
\fsquare(0.5cm,A) \otimes 
\fsquare(0.5cm,a_x) \otimes \cdots \otimes 
\fsquare(0.5cm,a_2) \otimes
\fsquare(0.5cm,b_m)
\otimes \Longsquare(0.5cm,0.8cm,b_{m-1})
\otimes \cdots\\ & \otimes&  \Longsquare(0.5cm,0.8cm,b_{y+1}) \otimes
\fsquare(0.5cm,B) \otimes \fsquare(0.5cm,b_y) 
\cdots \otimes \fsquare(0.5cm,{{b_2}}) \otimes
\fsquare(0.5cm,{{a_1}}) \otimes \fsquare(0.5cm,{b_1})
\otimes \fsquare(0.5cm,c_z) \otimes \cdots.
\end{eqnarray*}
Due to the conditions 
 $a_i<b_{k-1}\leq b_k$ and $a_i<B\leq b_{y+1}$ for 
$2\leq k<i\leq x$ , we can repeat the transformations above and get

\begin{eqnarray*}
ME(S)\douchic \cdots \otimes \fsquare(0.5cm,a_n) \otimes \cdots
& \otimes&
 \Longsquare(0.5cm,0.8cm,a_{x+1})  \otimes  
\fsquare(0.5cm,A) \otimes 
\fsquare(0.5cm,b_m) \otimes
\Longsquare(0.5cm,0.8cm,b_{m-1})
\otimes \cdots \otimes  \Longsquare(0.5cm,0.8cm,b_{y+1}) \otimes
\fsquare(0.5cm,B) \otimes \fsquare(0.5cm,b_y) \otimes \cdots\\
& \otimes & \Longsquare(0.5cm,0.8cm,b_{x+1}) \otimes
\fsquare(0.5cm,{a_x}) \otimes 
\fsquare(0.5cm,{b_x}) \otimes
\cdots \otimes 
\fsquare(0.5cm,{a_2}) \otimes
\fsquare(0.5cm,{b_2}) \otimes
\fsquare(0.5cm,{a_1}) \otimes \fsquare(0.5cm,{b_1})
\otimes \fsquare(0.5cm,c_z) \otimes \cdots.
\end{eqnarray*}
It follows from the conditions $A<b_i\leq b_{i+1}$ for 
$i=y+1,\cd,m$ and  $B\leq A< b_{y+1}$ that 

\begin{eqnarray}
ME(S)\douchic \cdots \otimes 
\fsquare(0.5cm,a_n) \otimes \cdots &\otimes&
\Longsquare(0.5cm,0.8cm,a_{x+1}) \otimes 
\fsquare(0.5cm,b_m) \otimes 
\Longsquare(0.5cm,0.8cm,b_{m-1})
\otimes \cdots \otimes \Longsquare(0.5cm,0.8cm,{b_{y+2}}) \otimes
\fsquare(0.5cm,{A}) \otimes \Longsquare(0.5cm,0.8cm,{b_{y+1}}) \otimes
\fsquare(0.5cm,B)
\otimes \fsquare(0.5cm,b_y) \otimes \cdots \nn \\
& \otimes & \Longsquare(0.5cm,0.8cm,b_{x+1}) \otimes
\fsquare(0.5cm,a_x) \otimes 
\fsquare(0.5cm,b_x) \otimes
\cdots \otimes 
\fsquare(0.5cm,a_1) \otimes \fsquare(0.5cm,b_1)
\otimes \fsquare(0.5cm,c_z) \otimes \cdots\nn \\
\douchic \cdots \otimes \fsquare(0.5cm,a_n) \otimes
\cdots & \otimes & \Longsquare(0.5cm,0.8cm,a_{x+1}) \otimes 
\fsquare(0.5cm,b_m) \otimes 
\Longsquare(0.5cm,0.8cm,b_{m-1})
\otimes \cdots \otimes \Longsquare(0.5cm,0.8cm,b_{y+2}) \otimes
\fsquare(0.5cm,{A}) \otimes \fsquare(0.5cm,{B}) \otimes
\Longsquare(0.5cm,0.8cm,{b_{y+1}})
\otimes \fsquare(0.5cm,b_y) \otimes \cdots \nn \\
& \otimes & \Longsquare(0.5cm,0.8cm,b_{x+1}) \otimes
\fsquare(0.5cm,a_x) \otimes 
\fsquare(0.5cm,b_x) \otimes
\cdots \otimes 
\fsquare(0.5cm,a_1) \otimes \fsquare(0.5cm,b_1)
\otimes \fsquare(0.5cm,c_z) \otimes \cdots 
\otimes \fsquare(0.5cm,c_1) \otimes \cdots 
\label{saigo}
\end{eqnarray}
Now, let us see the following Claim 1-3:

{\bf Claim 1.}
In (\ref{saigo}) one can find that $A$ and $B$ ($A\geq B$) are 
neighboring each other. Thus, 
applying the column bumping of (\ref{saigo}),
by the column bumping lemma (Lemma \ref{col-lem}) we obtain 

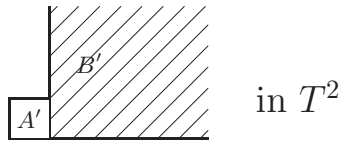
\setlength{\unitlength}{1pt}
\begin{picture}(90,90)
\put(10,10){\framebox(15,15){$A'$}}
\put(25,10){\line(1,0){60}}
\put(25,25){\line(0,1){35}}
\put(35,35){$B'$}
\hatchpattern{8}{\path}
    \hatch(25,10)(85,60)
\put(100,20){ \Large{in $T^2$}}

\end{picture}

where $A':=New(A)$ and $B':=New(B)$. 

{\bf Claim 2.}
Next, in the column bumping of $ME(S)$, since 
$a_1\leq \cd\leq a_x\leq A$, by the column bumping lemma 
(Lemma \ref{col-lem})
the new boxes by $a_1,\cd,a_x$ are placed 
on the right-side of $A'$. Similarly, since 
$c_1\leq \cd\leq c_z\leq b_1\leq \cd\leq b_y\leq B$, the new boxes 
by $c_1,\cd,c_z,b_1,\cd,b_x$ are placed on the 
right-side of $B'$.

{\bf Claim 3.}
As the definition of the map $\mathcal S_3$, 
the tableau $T^1$ is the
recording tableau of $T^2$. Then, it follows from Claim 2 that 
there are $x$ entry $k$'s on the right-side of $A'$ and 
$z+y$ entries $k+1$'s on the right-side of the same place 
as $B'$ in $T^1$.
We also know from Claim 1 that $B'$ is on the right-side of $A'$ 
and then there exist $z+y+1$ entry $k+1$'s on the right-side of
$A'$.

In $ME(Y_{\lm^1})\ot ME(T^1)$ let $n_1$ (resp. $n_2$) be the number of
$k$ (resp. $k+1$) on the left-side of $A'$.
Claim 3 implies that 
\begin{equation}
n_1=\lm^1+x,\qq
n_2=\lm^1+z+y+1.\label{n1n2}
\end{equation}
Since $z=\lm^1_k-\lm^1_{k+1}$ and $x\leq y$, one gets
\[
 n_2-n_1=(\lm^1_{k+1}+z+y+1)-(\lm^1_k+x)\geq 1,
\]
which contradicts that $T^1\in{\bf B}(\mu)_{\lm^1}^{\nu^1}$ and 
the case $S_{k,j}\geq S_{k+1,j}$ never occur. 
Thus, $S$ is a skew tableau.
It is immediate from the definition of ${\mathcal C}_2$
that $w(S)\,\,\douchik \,\,w(T)$, which means
$S$ is a Littlewood-Richardson skew tableau 
 and then $\mathcal C_2$ is well-defined. \qed

\nd
{\bf Proof of well-definedness of $\mathcal C$}
For the purpose we may show that $f$ is bijective, $f$ 
and $f^{-1}$ are
PJ-picture. The bijectivity of $f$ is obtained by the similar way
to that in \cite{NS,NS2}.
In order to show that $f$ 
and $f^{-1}$ are PJ-picture, we may see
for any $(i,j),(i,j+1), (i+1,j)\in \kappa^1$ and 
any $(a,b),(a,b+1), (a+1,b)\in\kappa^2$,
\begin{eqnarray*}
f(i,j)\leqslant_J f(i,j+1),\q
f(i,j)\leqslant_J f(i+1,j),\q
f^{-1}(a,b)\leqslant_J f^{-1}(a,b+1),\q
f^{-1}(a,b)\leqslant_J f^{-1}(a+1,b),
\end{eqnarray*}
These are also shown by the similar way to those in \cite{NS,NS2}.\qed

\section{Bijectivity of $\mathcal S$ and $\mathcal C$}

It suffices to show that $\mathcal C\circ\mathcal S={\rm id}$
and $\mathcal S\circ\mathcal C={\rm id}$.
To carry out this, we shall prove that 
$\mathcal C_i\circ\mathcal S_i={\rm id}$
and $\mathcal S_i\circ\mathcal C_i={\rm id}$ for $i=1,2,3$.
\subsection{$\mathcal S_1$ and $\mathcal C_1$}

Take $S\in {\mathbf S}(\kappa^1,\kappa^2)$ and  set 
$S':=\mathcal S_1\circ\mathcal C_1(S)$.
We have $\mathcal C_1(S)(i,j)=(S_{ij}, 
\lm^2_{S_{ij}}+p(S;i,j))$. Hence, by the definition of 
$\mathcal S_1$ we have $S'_{ij}=S_{ij}$, which implies $S'=S$ 
and then $\mathcal S_1\circ\mathcal C_1={\rm id}$.

For $f\in{\mathbf P}(\kappa^1,\kappa^2)$, set 
$g:=\mathcal C_1\circ\mathcal S_1(f)$.
The following lemma can proved similarly to 
\cite[Lemma 5.2]{NS} or \cite[Lemma 5.4]{NS2}.
\begin{lem}
Set $S=\mathcal S_1(f)$. Considering $Y_{\lm^2}\ot
ME(S)$, the entry $\fsquare(5mm,S_{ij})$ is added to
the position $f(i,j)\in\kappa^2$.
\end{lem}

Since $S_{ij}=f_1(i,j)$ and 
$g(i,j)=(S_{ij},\lm^2_{S_{ij}}+p(S;i,j))$, we get 
$g_1(i,j)=f_1(i,j)$.
We know that $S_{ij}(=k)$ is the $p(S;,i,j)$-th 
entry equal to $k$ and $f_2(i,j)
=\lm^2_{S_{ij}}+p(S;,i,j)=g_2(i,j)$, which shows $f=g$ and then
$\mathcal C_1\circ\mathcal S_1={\rm id}$.

\subsection{$\mathcal S_2$ and $\mathcal C_2$}

Set $w':=\mathcal S_2\circ\mathcal C_2(w)$ for 
$w\in{\mathbf W}(\kappa^1,\kappa^2)$ and write
\[
 w=\begin{pmatrix}w_1\\ w_2\end{pmatrix}
=\begin{pmatrix}
d_1d_2\cd d_N\\ c_1c_2\cd c_N\end{pmatrix},
\qq\qq
w'=\begin{pmatrix}w'_1\\ w'_2\end{pmatrix}=\begin{pmatrix}
b_1b_2\cd b_N\\ a_1a_2\cd a_N\end{pmatrix}.
\]
Note that number of $i$ in $w_1$ is just equal to $\kappa^1_i$.
For $S:=\mathcal C_2(w)$, we have $ME(S)=
\fsq(c_1) \otimes \fsq(c_2) \otimes \cdots \otimes \fsq(c_n)$
and then $w_2=w'_2$ by the definition of $\mathcal S_2$.
The number $b_i$ is the row number of $a_i$ in $S$.
Thus, since the number of $i$ in $w'_1$ is $\kappa^1_i$, 
$d_1\leq \cd\leq d_N$ and $b_1\leq \cd\leq b_N$, 
we have $w_1=w'=1$ and then $w=w'$, which means 
$\mathcal S_2\circ \mathcal C_2={\rm id}$.

It is trivial from the definition of the maps
$\mathcal S_2$ and $\mathcal C_2$ that
$\mathcal C_2\circ\mathcal S_2={\rm id}$.

\subsection{$\mathcal S_3$ and $\mathcal C_3$}

We have seen the well-definedness of the maps
$\mathcal S_3$ and $\mathcal C_3$ and these maps 
are certain restriction of usual RSK correspondence 
in terms of column bumping. Thus, we obtain 
$\mathcal S_3\circ\mathcal C_3={\rm id}$ and 
$\mathcal C_3\circ\mathcal S_3={\rm id}$.

Now, we obtain 
$\mathcal S_i\circ\mathcal C_i={\rm id}$ and 
$\mathcal C_i\circ\mathcal S_i={\rm id}$ ($i=1,2,3$) and then 
$\mathcal S\circ\mathcal C={\rm id}$ and 
$\mathcal C\circ\mathcal S={\rm id}$.
So, we have completed the proof of Theorem \ref{main}.
\qed

\end{document}